\documentclass[12pt,twoside, epsf]{article}

\usepackage{color,graphicx,times}

\usepackage{amssymb}
\usepackage{amsmath}
\usepackage{amscd}
\usepackage{float}
\usepackage{graphicx}
\usepackage{amsfonts}
\usepackage{pb-diagram}
\usepackage{eufrak}

 \newcommand{\CrossGlyph}{\raisebox{-0.25\height}{\includegraphics[width=0.5cm]{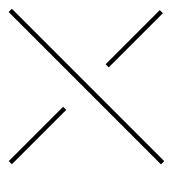}}}
\newcommand{\HSmoothGlyph}{\raisebox{-0.25\height}{\includegraphics[width=0.5cm]{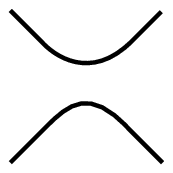}}}
\newcommand{\VSmoothGlyph}{\raisebox{-0.25\height}{\includegraphics[width=0.5cm]{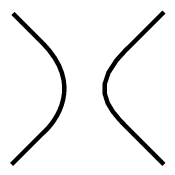}}}
\newcommand{\VirtualGlyph}{\raisebox{-0.25\height}{\includegraphics[width=0.5cm]{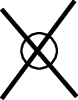}}}

\newcommand{\IMGlyph}{\raisebox{-0.25\height}{\includegraphics[width=0.5cm]{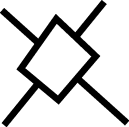}}}
\newcommand{\IMGlyphL}{\raisebox{-0.25\height}{\includegraphics[width=0.5cm]{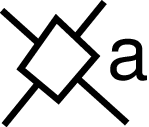}}}

\newcommand{\CDotDiag}{\raisebox{-0.25\height}{\includegraphics[width=0.5cm]{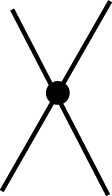}}}
\newcommand{\CCircleDiag}{\raisebox{-0.25\height}{\includegraphics[width=0.5cm]{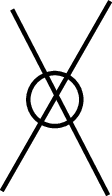}}}
\newcommand{\CBoxDiag}{\raisebox{-0.25\height}{\includegraphics[width=0.5cm]{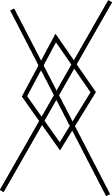}}}
\newcommand{\CCircleDiagL}{\raisebox{-0.25\height}{\includegraphics[width=0.5cm]{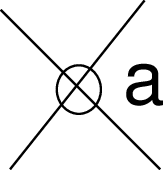}}}

\makeatother

\let \ttorg \tt \def \tt{\ttorg \obeyspaces}

\begin{document}

\date{}

\title{\bf  Strict Equivalence of Multi-Virtual Linkoids }

\author{Louis H. Kauffman \\ 
Department of Mathematics, Statistics and Computer Science \\ 
University of Illinois at Chicago \\ 
851 South Morgan Street\\
Chicago, IL, 60607-7045\\
and\\
International Institute for Sustainability with Knotted Chiral Meta Matter\\
$WPI-SKCM^{2}$\\
Hiroshima University, 1-3-1 Kagamiyama\\
Higashi-Hiroshima, Hiroshima 739-8531\\
Japan}

\maketitle

\thispagestyle{empty}

\begin{abstract}
We utilize  multi-virtual knot theory where there are a multiplicity of virtual crossings to study {\it strict virtual linkoids}. In strict virtual linkoid theory, local moves  define all virtual moves and Reidemeister moves. 
In the strict equivalence, no moves, classical or virtual, can transfer an arc across a linkoid endpoint. By taking closures of strict virtual linkoids that are multi-virtual knots and links, we obtain new invariants for strict virtual linkoids.
Generalized bracket polynomial invariants and generalized loop bracket polynomial invariants (for planar strict virtual linkoids) are studied in this context.
The paper defines {\it virtual  polar links} where there are degree two nodes in virtual link diagrams across which isotopies are forbidden. The paper shows how multi-virtual theory and its concepts 
can be applied to obtain invariants for polar virtual links.
 \end{abstract}

\noindent {\bf MSC2020-Mathematics Subject Classification System:} $05C15, 57K10, 57K12, 57K14.$\\

\section{Introduction}
This paper initiates a study of strict virtual linkoid theory, where local moves define all virtual moves and Reidemeister moves. 
In the strict equivalence, no moves, classical or virtual, can transfer an arc across a linkoid endpoint. By taking closures of strict virtual linkoids as multi-virtual knots and links, we obtain new invariants for strict virtual linkoids.
The paper defines strict polar multi-virtual links where there are degree two nodes in the diagrams across which isotopies are forbidden. Generalized bracket polynomial invariants and generalized loop bracket polynomial invariants (for planar strict virtual linkoids) are studied in this context. The paper shows how multi-virtual theory and its concepts 
can be applied to obtain invariants for polar virtual links. \\

This paper uses a generalization of virtual knot theory that we call
 multi-virtual knot theory \cite{MVKT,MVALG}. Multi-virtual knot theory invokes a multiplicity of types of virtual crossings. We begin the paper with a review of multi-virtual knot theory and the generalized bracket polynomials that occur in this theory.
 We then introduce strict virtual knotoids and linkoids and show how they can be studied using multi-virtual knot theory. We end the paper with a generalization of polar knots and links to strict virtual polar knots and links and how the methods of 
 this paper can be used to study virtual polar links.\\ 
 
Section 2 begins the definition and exploration of {\it multi-virtual knot theory.} These virtual crossings all detour over one another, while the classical crossings do not detour over the virtual crossings. The reader may wish to read  earlier introductions to virtual knot theory, multi-virtual knot theory and related research by the author \cite{VKT,DVK,vkt,rotvkt,MVKT,MVALG}. Multi-virtual knot theory is presented here in terms of knot and link diagrams that have
classical (and sometimes flat classical) crossings and virtual crossings. The classical crossings interact with one another via the Reidemeister moves and the virtual crossings interact via detour moves. In a detour move a consecutive sequence of virtual crossings of the same type can be excised from the diagram and the endpoints of this excision can be reconnected by any diagrammatic arc that intersects the rest of the diagram in a sequence of virtual crossings of this same type. In sections 2.1, 2.2, 2.3 and 2.4 we discuss this definition of multi-virtual knot theory and we give a sketch of the topological interpretation of this definition in the case of single virtual knot theory.
Single virtual knot theory can be interpreted in terms of knots and links embedded in thickened surfaces taken up to knot theoretic equivalence in these surfaces and up to one-handle stabilization. In Sections 2.3 and 2.4 we describe a similar but more restricted topological interpretation for multi-virtual knot theory that is obtained by adding labeled handles at the site of each virtual crossing. Handles with the same label can merge, and handle stabilization arises in relation to these operations.\\

In {\it strict multi-virtual knot theory} the detour moves are replaced by local detour moves as shown in Figure~\ref{Figure 1}. In the case of closed loop diagrams the strict theory is the same as the theory defined by general detour moves. That is, for closed loop diagrams every detour move can be written as a composition of local detours. But when we work with diagrams that have endpoints as with tangles or knotoids or linkoids, then the restriction to local moves makes a difference that leads to the considerations of this paper.\\

Section 2.5 defines a generalized bracket polynomial for multi-virtual knot theory. The basic bracket expansion is given by the formula
$$\langle \CrossGlyph \rangle = A \langle \HSmoothGlyph \rangle + B \langle \VSmoothGlyph \rangle.$$
See \cite{LKNew, LKStat,LKTutte,KauffJaeger,Kauff:KP}. By itself, the bracket polynomial yields a significant connection between low dimensional topology and the graph theoretical structures of the chromatic, dichromatic and Tutte polynomials. Here we extend the bracket to an invariant of multi-virtual knots. The most general version uses $B=A^{-1}$ and $\delta = -A^{2} - A^{-2}$ to give an invariant of regular isotopy of multi-virtual knots and links. In this formulation we have a state summation formula  in the form $$\langle K \rangle = \sum_{S} \langle K | S \rangle \langle S \rangle$$ where $ \langle K | S \rangle$ denotes a product of $A's$ amd $B's$ corresponding to a choice of smoothing of each crossing in the diagram $K.$ The state $S$ indicates such a choice of smoothings. $S$ consists in a collection of curves in the plane intersecting transversely in virtual crossings of different types. $\langle S \rangle$ denotes the equivalence class of $S$ under detour moves on these virtual crossings. In general these equivalence classes can be difficult to determine, but in many special cases the generalized invariant can yield much information.\\

We study a special case of this invariant, the {\it chromatic bracket}, by combining it with a graphical deletion contraction algorithm that is applicable to multi-virtual diagrams with two distinct virtual crossings. When there are no more than two virtual crossings we denote them by round glyphs and square glyphs. The equations are as follows.
$$\langle \CrossGlyph \rangle = A \langle \HSmoothGlyph \rangle + B \langle \VSmoothGlyph \rangle$$
$$\langle O \rangle = \delta$$
$$\langle \IMGlyph \rangle  = 2 \langle \CDotDiag \rangle  - \langle \VirtualGlyph \rangle $$
In the last equation, a square virtual glyph is expanded as twice a black node minus a round virtual glyph. The dark node is the analog of graphical contraction and the round node is the analog of graphical deletion. Any collection of state loops that are connected by black nodes is evaluated as $\delta.$ Loops with common round nodes are regarded as disjoint so that $n$ loops with no black nodes will evaluate as $\delta^{n}.$ We recall in detail (below) the proof from \cite{MVKT} that the chromatic bracket polynomial is a multi-virtual link regular isotopy invariant for two virtual crossing types when $B=A^{-1}$ and $\delta = -A^{2} - A^{-2}.$ 
Section 2.5 gives many topological examples for these polynomials. \\

Section 3 introduces the concept of {\it strict} knotoids and linkoids, where the virtual detour moves are restricted to compositions of local moves that never pass an arc across an endpoint of the knotoid or linkoid. This  restricted
equivalence for knotoids and linkoids has a rich mathematical structure as we show in the rest of the paper. This section shows how standard classical under and over closures for linkoids and knotoids are not well defined in the virtual category, but that a special form of virtual closure ({\it box closure}) is effective where we add a new virtual crossing in performing the closure. One can then compute generalized brackets for multi-virtual knots for these closures to obtain new invariants of strict virtual knotoids and linkoids.
A number of examples are given, and we use known limitations of the Jones polynomial to give examples of strict linkoids whose nontriviality cannot be detected by the generalized bracket polynomials. The section ends by showing how the box closure can be used to estimate the complexity (height) of a strict virtual knotoid.\\

Section 4 discusses planar strict linkoids. Here the linkoids are diagrammed in a plane and so arcs cannot be swung around the two-sphere. We explain how the generalized bracket polynomial becomes a generalized loop bracket \cite{Neslihan,Turaev} with extra combinatorial states due to the strict detour equivalence. In the original loop bracket, new variables appear from irreducible states corresponding to nested loops surrounding a knotoid arc. In the generalized loop bracket there are infinitely many new combinatorial states. See Figure~\ref{examples} and the discussion in this section. Classification of planar multi-virtual states for linkoids will be the subject of a sequel to the present paper.\\

Section 5 introduces the concept of a virtual polar knot or link, where there are degree two nodes in the diagram across which no arcs can move in an isotopy. It is shown how the techniques of this paper apply to obtaining invariants of virtual polar knotoids
and linkoids.\\

Section 6 discusses problems and prospects that arise from the present work.\\

\noindent {\bf Acknowledgement.} The author thanks the Mathematisches  Forschungsinstitut Oberwohlfach for its hospitality during crucial stages in the preparation of this paper.\\

\section{Virtual and multi-virtual Knot Theory}
In this section we will give the basics of a generalization of virtual knot theory to {\it multi-virtual knot theory} where threre is a multiplicity of types of virtual crossing. This is further generalized to knotoids and linkoids.
We begin with a description of standard virtual knot theory and how it will be modified.  There are many variants of knot theory that occur in relation to virtual knot theory. The reader can see a survey of these variants in \cite{VKT,vkt,SL,DVK,DKT,DKK,rotvkt,affineindex,VKC,VKC1,SL,IP,VSA,FennBook,KFM}. The paper by the author on multi-virtual knot Theory \cite{MVKT} is a direct source of background material for this paper and the joint work \cite{MVALG} gives further information about algebraic invariants of multi-virtual knots.In this section we repeat material from \cite{MVKT} in order to obtain a smooth exposition of the basics.\\

\subsection{Virtual Knot Theory}
Knot theory
studies the embeddings of curves in three-dimensional space.  Standard virtual knot theory studies the  embeddings of curves in thickened surfaces of arbitrary
genus, up to 1-handle stabilization. Virtual knots have a special diagrammatic theory, described below,
that makes handling them
very similar to the handling of classical knot diagrams.  \bigbreak  

\begin{figure}[htb]
     \begin{center}
     \begin{tabular}{c}
     \includegraphics[width=10cm]{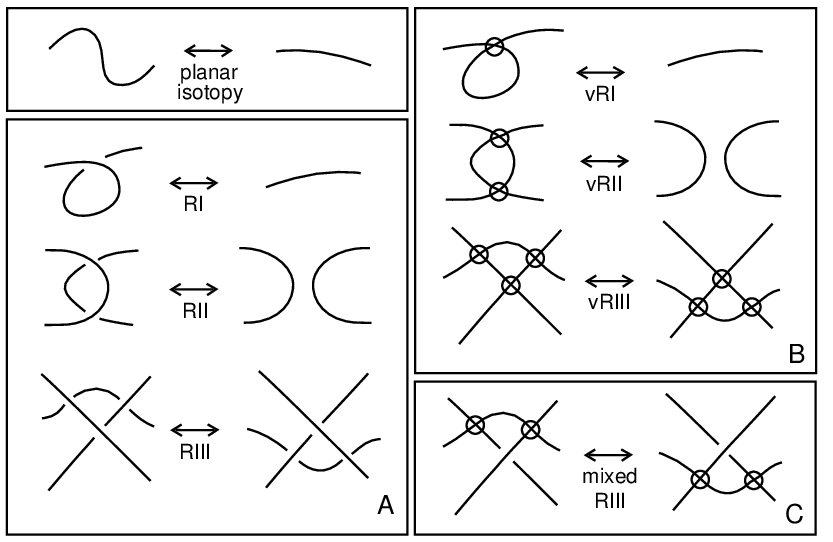}
     \end{tabular}
     \caption{\bf Moves}
     \label{Figure 1}
\end{center}
\end{figure}

In the diagrammatic theory of virtual knots one adds 
a single type of {\em virtual crossing} (see Figure~\ref{Figure 1}) that is neither an over-crossing
nor an under-crossing.  A virtual crossing is represented by two transversely crossing segments, with a small circle
placed around the crossing point. 
\bigbreak

Moves on virtual diagrams generalize the Reidemeister moves for classical knot and link
diagrams.  See Figure~\ref{Figure 1}.  One can summarize the moves on virtual diagrams by saying that the classical crossings interact with
one another according to the usual Reidemeister moves while virtual crossings are artifacts of the attempt to draw the virtual structure in the plane. 
A segment of diagram consisting of a sequence of consecutive virtual crossings can be excised and a new connection made between the resulting
free ends. If the new connecting segment intersects the remaining diagram (transversely) then each new intersection is taken to be virtual.
Such an excision and reconnection is called a {\it detour move}.
Adding the global detour move to the Reidemeister moves completes the description of moves on virtual diagrams. In Figure~\ref{Figure 1} we illustrate a set of local
moves involving virtual crossings. The global detour move is
a consequence of  moves of type (B) and (C) in Figure~\ref{Figure 1}. The detour move is illustrated in Figure~\ref{Figure 2}.  Virtual knot and link diagrams that can be connected by a finite 
sequence of these moves are said to be {\it equivalent} or {\it virtually isotopic}.
\bigbreak

\begin{figure}
     \begin{center}
     \begin{tabular}{c}
     \includegraphics[width=8cm]{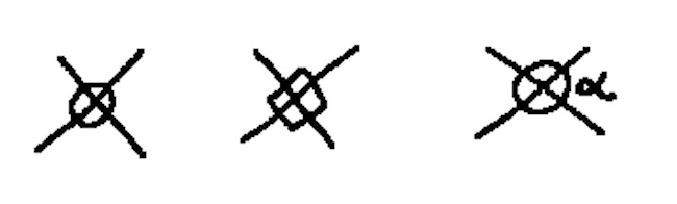}
     \end{tabular}
     \caption{\bf Virtual Crossing Notations}
     \label{EFF1}
\end{center}
\end{figure}

\begin{figure}
     \begin{center}
     \begin{tabular}{c}
     \includegraphics[width=8cm]{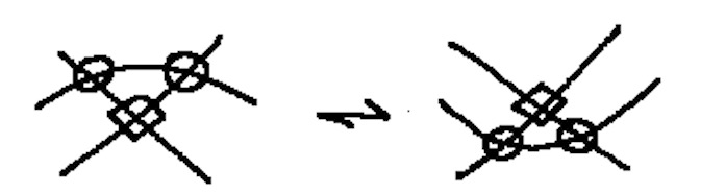}
     \end{tabular}
     \caption{\bf Local Detour Move for Virtuals across a Virtual}
     \label{EFF3}
\end{center}
\end{figure}

\begin{figure}[htb]
     \begin{center}
     \begin{tabular}{c}
     \includegraphics[width=10cm]{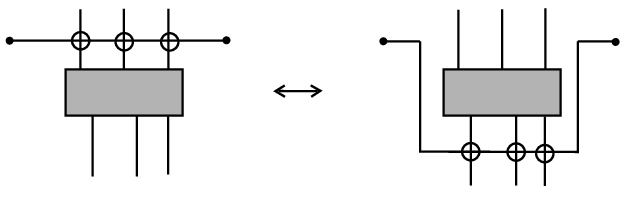}
     \end{tabular}
     \caption{\bf Detour Move}
     \label{Figure 2}
\end{center}
\end{figure}

\begin{figure}[htb]
     \begin{center}
     \begin{tabular}{c}
     \includegraphics[width=10cm]{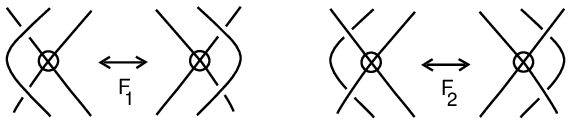}
     \end{tabular}
     \caption{\bf Forbidden Moves}
     \label{Figure 3}
\end{center}
\end{figure}

Another way to understand virtual diagrams is to regard them as representatives for oriented Gauss codes \cite{GPV}, \cite{VKT,DVK} 
(Gauss diagrams). Such codes do not always have planar realizations. An attempt to embed such a code in the plane
leads to the production of the virtual crossings. The detour move makes the particular choice of virtual crossings 
irrelevant. {\it Virtual isotopy is the same as the equivalence relation generated on the collection
of oriented Gauss codes by abstract Reidemeister moves on these codes.}  
\bigbreak

Figure~\ref{Figure 3} illustrates  two {\it forbidden moves}. Neither of these follows from Reidmeister moves plus detour move, and 
indeed it is not hard to construct examples of virtual knots that are non-trivial, but will become unknotted on the application of 
one or both of the forbidden moves. The forbidden moves change the structure of the Gauss code and, if desired, must be 
considered separately from the virtual knot theory proper. \\

Figure~\ref{EFF1} and Figure~\ref{EFF3} illustrate how we may have more than one type of virtual crossing and that when we have such a multiplicity of virtual crossings we allow (in the framework of the present work) that distinct virtual types can detour over one another.\\

Later in the paper we will restrict virtual moves and detour moves to local Reidemeister-type moves such as shown in the right of Figure~\ref{Figure 1} and in Figure~\ref{EFF3}. When the virtual moves are restricted to compositions of local moves, the equivalence relation is called {\it strict virtual isotopy}. As we shall see, strict virtual isotopy for knotoids and linkoids (described in section 3 and beyond) gives a new domain for exploration.\\
\bigbreak

\subsection{Interpretation of Virtuals Links as Stable Classes of Links in  Thickened Surfaces}
There is a topological interpretation \cite{VKT,DVK,Carter,Kamada3} for virtual theory in terms of embeddings of links
in thickened surfaces.  Regard each 
virtual crossing as a shorthand for a detour of one of the arcs in the crossing through a 1-handle
that has been attached to the 2-sphere of the original diagram.  Interpret each virtual crossing to obtain an embedding of a collection of circles into a thickened surface  $S_{g} \times R$ where $g$ is the 
number of virtual crossings in the original diagram $L$, $S_{g}$ is a compact oriented surface of genus $g$
and $R$ denotes the real line.  Two such surface embeddings are called
{\em stably equivalent} if one can be obtained from another by isotopy in the thickened surfaces, 
homeomorphisms of the surfaces and the addition or subtraction of empty handles (the knot does not go through the handle).

\begin{figure}
     \begin{center}
     \begin{tabular}{c}
     \includegraphics[width=10cm]{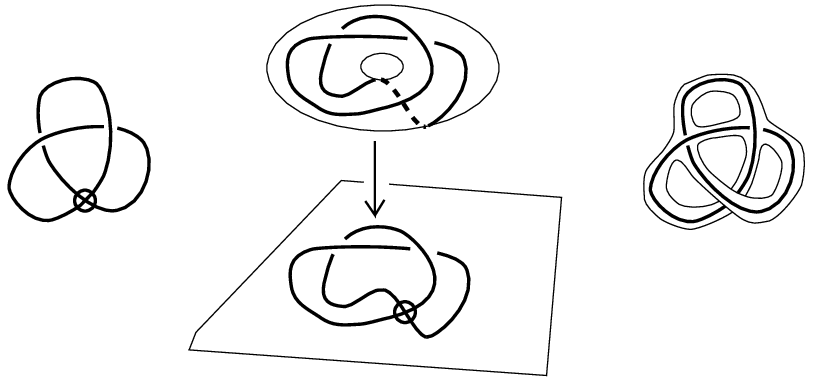}
     \end{tabular}
     \caption{\bf Surfaces and Virtuals}
     \label{Figure 4}
\end{center}
\end{figure}

\noindent We have the
\smallbreak
\noindent
{\bf Theorem 1 \cite{VKT,DKT,DVK,Carter}.} {\em Two virtual link diagrams are isotopic if and only if their corresponding 
surface embeddings are stably equivalent.}  
\smallbreak
\noindent
\bigbreak  

\noindent In Figure~\ref{Figure 4} we illustrate this association of virtual diagrams and knot and link diagrams on surfaces.
The projection of the knot diagram on the torus to a diagram in the plane (in the center of the figure) has a virtual crossing in the 
planar diagram where two arcs (that do not form a weave on the surface) project to the same point in the plane. Virtual 
crossings are seen as artifacts of projection. The figure shows a virtual diagram on the left and an ``abstract knot diagram" \cite{Kamada3,Carter} on the right.
The abstract knot diagram is a realization of the knot on the left in a thickened surface with boundary. The abstract knot diagram is obtained by making a neighborhood of the 
virtual diagram that resolves the virtual crossing into arcs that travel on separate bands. The virtual crossing appears as an artifact of the
projection of the surface to the plane. See \cite{VKT,DKT} and other papers by the author and in
the literature of virtual knot theory.
\bigbreak
 
\subsection{Multi-Virtual Knot Theory}
We use the following notations:
ordinary virtual crossings are denoted by $\CCircleDiag,$ and a second type of virtual crossing is
denoted by $\CBoxDiag.$  
In principle, we can have as many distinct types of virtual crossings as we like and they can be labeled by glyphs like the above with 
a letter attached and the letter running over a chosen index set as in $\CCircleDiagL,$ or $\IMGlyphL.$\\

\noindent {\bf Detour Axiom for multi-virtual Crossings.} Each virtual crossing type detours over all the other types of virtual crossings and over the classical crossings.
Figures  \ref{Figure 2}, \ref{EFF3} and  \ref{EFF4} illustrate the detour moves.\\ 

\noindent {\bf Definition.} {\it Multi-Virtual Knot Theory (MVKT)} is the diagrammatic theory generated by diagrams with multi-virtual crossings, using the Detour Axiom stated above in conjunction with the 
classical Reidemeister moves. For a given instantiation of MVKT, a set of virtual crossing types is chosen with a method for indicating them. Note that in one application of the detour move, only one virtual crossing type will be used.
Thus if one excises a sequence of virtual crossings of a given type, then then the new connection will create only virtual crossings of that type. If the excision does not remove any virtual crossings, then the reconnection can use
any single virtual crossing type that is available.\\

Note the two component ``link" with two virtual crossings in Figure~\ref{EFF2}. This link is not equivalent to a disjoint union of two circles. We will see a proof of this below. The figure illustrates the generalized bracket expansion of a multi-virtual knot with one classical crossing and two virtual crossings. We will define the generalized bracket below. \\

Figure~\ref{onecomp} shows a two-virtual diagram of one component and a series of detour moves that reduce it to a circle.\\

 \noindent {\bf Conjecture.} Any single component diagram (as in Figure~\ref{onecomp}), decorated with no more than two distinct virtual crossings, reduces by detour moves to a circle.\\

The author has previously considered multi-virtual crossings with hierarchies of detour moves among them and there is a paper \cite{MFlat} taking steps with this idea.
In multi-virtual theory, all the virtual crossings are on the same footing. They can detour over one another, but classical crossings cannot detour over virtual crossings.\\

 Many questions arise about the relationship of classical knot theory, virtual knot theory and multi-virtual knot theory. It is known that classical knot theory embeds in virtual knot theory in the sense that if two classcial knot or link diagrams are equivalent as virtual diagrams, then they are equivalent as classical diagrams \cite{VKT}. The same result holds for multi-virtual knot theory.\\
 
 \noindent {\bf Theorem.} Let $K$ and $L$ be classical link diagrams in classical knot theory (KT), and let MVKT denote a multi-virtual knot theory (with a given cardinality of virtual crossing types). Since $K$ and $L$ can be viewed as diagrams in MVKT we have a tautological mapping $F: KT \longrightarrow MVKT $ where F(K) is the same diagram seen in the larger category of virtual moves. Then if $F(K)$ and $F(L)$ are equivalent in MVKT, they are equivalent in KT. Thus classical knot theory embeds in multi-virtual knot theory.\\

\noindent {\bf Proof.} The proof has the same structure as the proof for virtual knot theory in \cite{VKT}. It is known from the work of Waldhausen that the fundamental group and periperal subgroup as a pair $(G,P)$ classify classical knots and links. Thus if $K$ and $L$ are classical links with peripheral subgroup pairs $(G(K), P(K))$ and  $(G(L), P(L))$ then $K$ is classically ambient istopic to $L$ if and only if the corresponding group pairs are isomorphic. The group pair of a classical link or a virtual link diagram $K$ can be given by a Wirtinger presentation read from the diagram in conjunction with longitudes $\lambda_{i}(K)$ (one for each link component) that is obtained by walking along a component and writing the product of the Wirtinger generators 
(corrsponding to oriented arcs of the diagram) that are underpassed in the course of the walk. These longitudes are not affected by detour moves. From this it follows that if $L$ is obtained from $K$ by a sequence of Reidemeister moves and virtual detour moves (in the sense of MVKT), then  $(G(L), P(L))$  is isomorphic as a pair with $(G(K), P(K)).$ If $K$ and $L$ are classical diagrams, it follows that they are classically equivalent. This completes the proof of the Theorem.  $\hfill\Box$

\begin{figure}
     \begin{center}
     \begin{tabular}{c}
     \includegraphics[width=10cm]{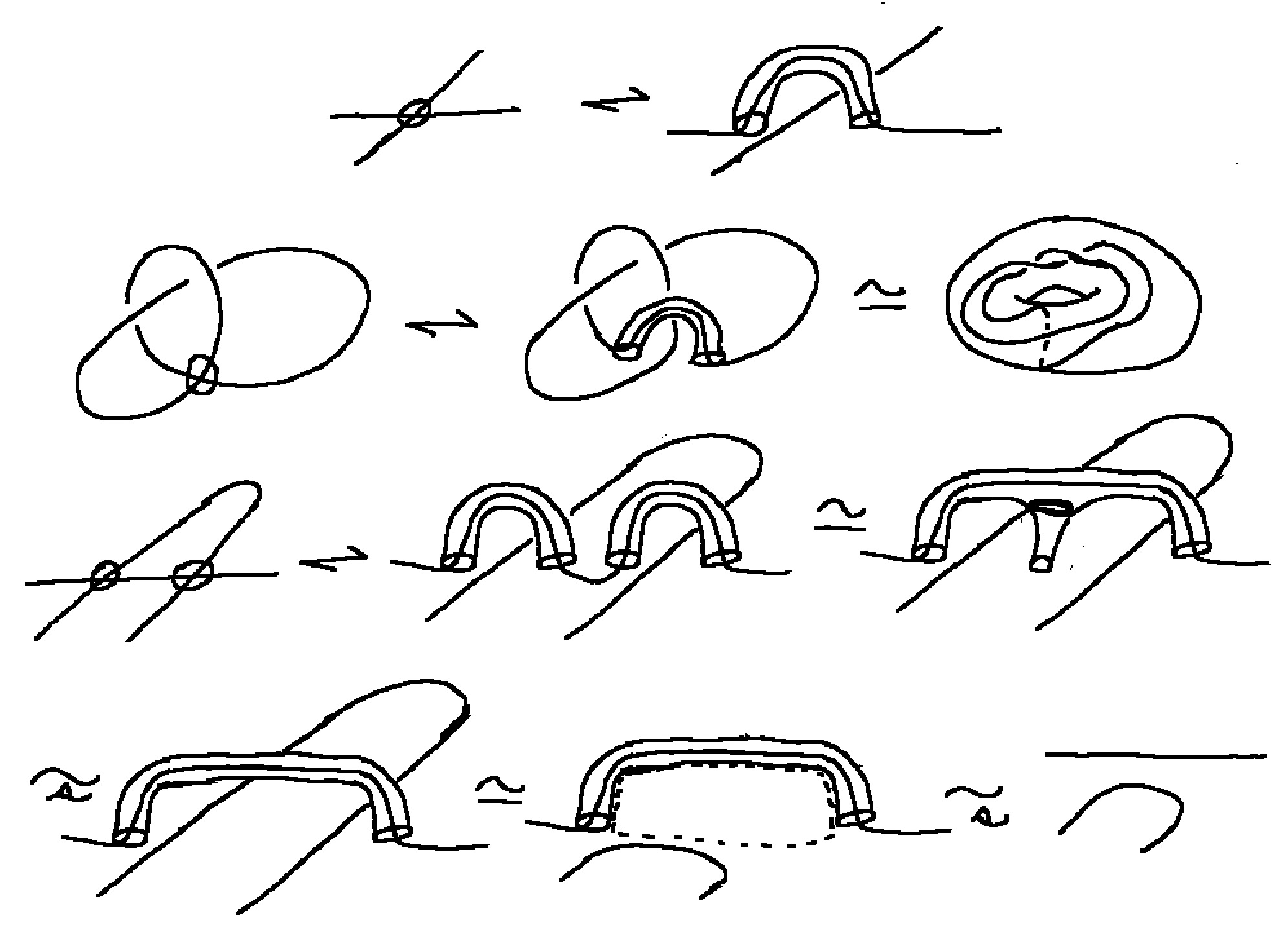}
     \end{tabular}
     \caption{\bf Replacing Virtual Crossings by Handle Detours}
     \label{Handles}
\end{center}
\end{figure}

\subsection{Surface Interpretations for multi-virtual Knot Theory}
Single type virtual knot theory can be interpreted in terms of knots in thickened surfaces. One obtains a knot or link in a surface by using an abstract link diagram, obtained by producing a ribbon neighborhood of the virtual diagram and using two local ribbons for each virtual crossing as shown in Figure~\ref{Figure 4}. Another method that associates a surface to a virtual diagram is illustrated in Figure~\ref{Handles}. In that figure we add a handle to the two-sphere or plane for the virtual link at each virtual crossing. The virtual crossing is replaced by two arcs, one going through the handle, and one going underneath the handle. This construction can be generalized to multi-virtual knot theory by labelling the resulting handles according to the label of the virtual crossing. The handles then interact as shown in Figure ~\ref{Handles} only when they have the same label. In the figure we have illustrated how a detour move (virtual Reidemeister move two) is seen in terms of handle stabilization. The two handles interact and undergo a 1-handle stabilization that fuses them into a single handle. This allows the arc passage for the move and then we destabilize the resulting handle, taking us back to the two sphere or the plane. Handles that interact come from the same type of virtual crossing.  This is a mixed combinatorial and topological interpretation.\\

\subsection{Generalized Bracket Polynomial}
We now start the proper discussion of multi-virtual knots and links. In Figure~\ref{EFF1} we illustrate graphical notations for virtual crossings.
The traditional circle notation is shown along with a local square notation for a distinct virtual crossing, and a circle notation with an index $\alpha.$
The $\alpha$ can run over any convenient index set, allowing any chosen cardinality of virtual crossings. \\

\begin{figure}
     \begin{center}
     \begin{tabular}{c}
     \includegraphics[width=10cm]{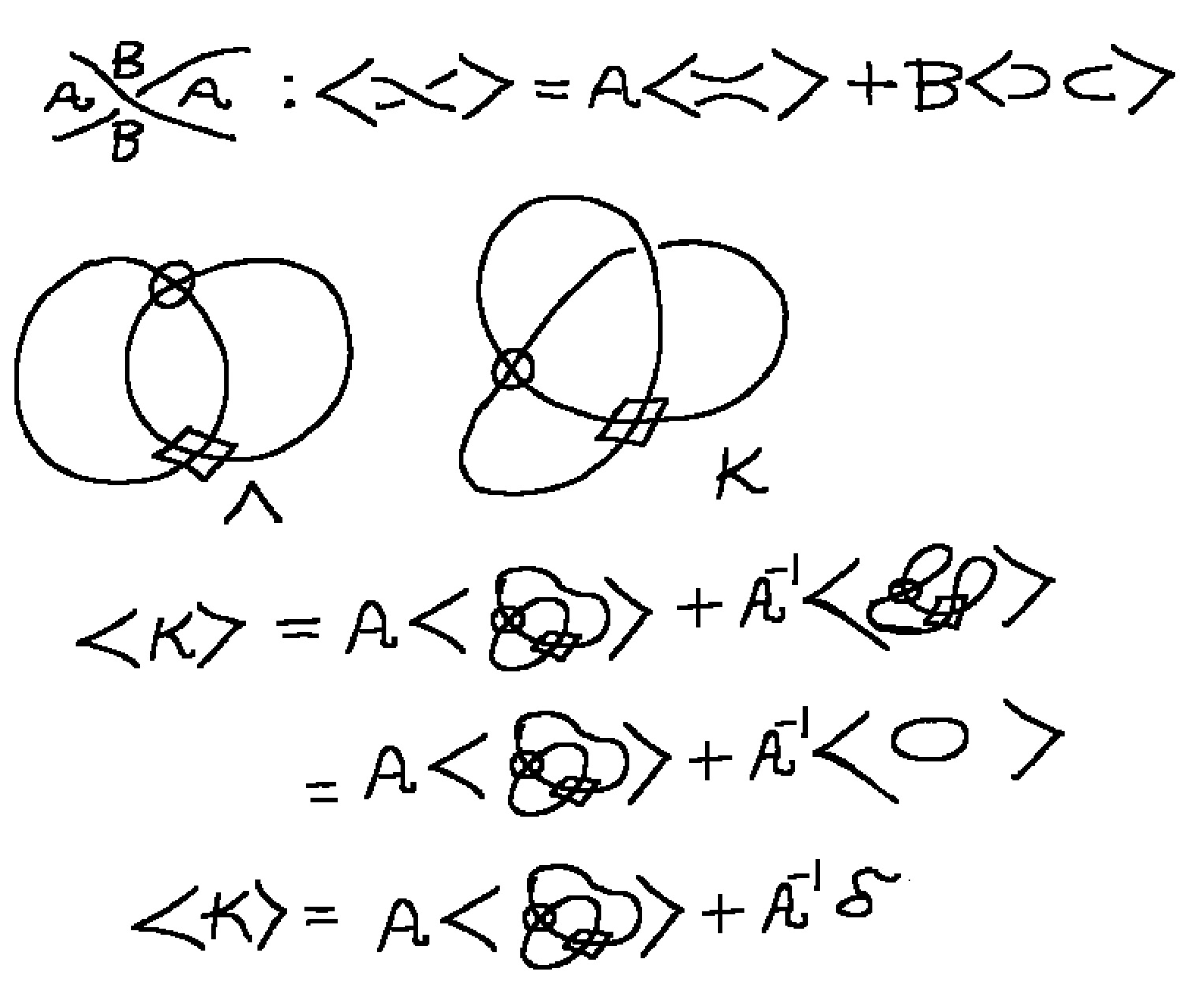}
     \end{tabular}
     \caption{\bf Double Virtual Link and Double Virtual Knot}
     \label{EFF2}
\end{center}
\end{figure}

In Figure~\ref{EFF2} we show the simplest multi-virtual link, a link with two virtual crossings, one of circular type and one of square type.
Our axioms for virtual crossings do not permit two different types to cancel, and  this link can be proved to be non-trivial. We will shortly explain an algebraic invariant that verifies this fact.\\

In Figure~\ref{EFF3} and Figure~\ref{EFF4} we illustrate a simple example of a detour move of two virtuals of one type, across a virtual of another type. In Figure~\ref{EFF3}  we show a detour move in relation to another virtual crossing. In Figure~\ref{EFF4}  we show a detour move in relation to a classical crossing and in relation to a 4-valent graphical node. We will use such nodes in our constructions for invariants.\\

Figure~\ref{EFF5} illustrates the non-cancellation of distinct virtual crossings.\\

We first define a generalized bracket polynomial for any multi-virtual theory. This bracket invariant has the usual skein relation and loop evaluation \cite{K,LKNew,LKStat,Kauff:KP}.
$$\langle \CrossGlyph \rangle = A \langle \HSmoothGlyph \rangle + A^{-1} \langle \VSmoothGlyph \rangle$$
$$\langle O \rangle = -A^2 -A^{-2}= \delta$$\\

Applying the expansion formula recursively will yield collections of loops with multi-virtual crossings. The generalized bracket has state expansion 
$$<K> = \Sigma_{S} A^{a(S)}B^{b(S)} \langle S \rangle$$ 
where $S$ is a state of the diagram $K$ obtained by choosing a smoothing for each 
crossing of $K$ and $a(S)$ is the number of smoothings of type $A,$ while $b(S)$ is the number of smoothings of type $B.$ Smoothing types correspond to the convention shown in 
Figure~\ref{EFF2} where the local regions of a crossing are labeled $A$ and $B$ so that the $A$ regions are swept by a counterclockwise turn of the overcrossing arc. A smoothing is of type $A$ if the two $A$ regions are joined by the smoothing. Each state $S$ in the virtual bracket expansion is a diagram with multi-virtual crossings. For the generalized bracket, we take $\langle S \rangle$ to be the multiple flat virtual class (i.e. the equivalence class of the state up to virtual detour moves)  of the state $S$ with the caveat that any free disjoint circle in $S$ (or in a diagram equivalent to $S$) is evaluated as $\delta.$ Thus we have that $\langle S O \rangle = \delta \langle S \rangle.$\\

Since, in general, we do not know full information about the class of a state $S$ up to detour moves, this generalized bracket is a {\it picture valued} invariant where the pictures (diagrams) of the states represent their 
equivalence classes. There are different ways to create invariants for these states, and so one leaves the generalized bracket in this open-ended form.\\

In Figure~\ref{EFF2} we illustrate the bracket expansion of a multi-virtual knot $K.$ In this expansion we obtain one picture $\langle \Lambda \rangle$ where $\Lambda$ is the state indicated in the figure, consisting of two circles intersecting in two distinct virtual crossings.\\

Determining the virtual detour class of these states in the multi-virtual setting is an additional mathematical problem. We make a special evaluation in the case of {\it double multi-virtual theory} where we use two virtual crossing types, one indicated by a circle and the other by a square as in Figure~\ref{EFF1}.
Then we take $\langle S \rangle$ to be the graphical evaluation to be discussed below using the expansion in Figure~\ref{EFF7} and Figure~\ref{EFF8}.
We call this specialized bracket the {\it chromatic bracket} since its extra rule is motivated by the coloring problems for graphs. The chromatic bracket is then governed by the rules shown below.\\

$$\langle \CrossGlyph \rangle = A \langle \HSmoothGlyph \rangle + B \langle \VSmoothGlyph \rangle$$
$$\langle O \rangle = \delta$$
$$\langle \IMGlyph \rangle  = 2 \langle \CDotDiag \rangle  - \langle \VirtualGlyph \rangle $$\\

In the last equation, a square virtual glyph is expanded as twice a black node minus a round virtual glyph. The dark node is the analog of graphical contraction and the round node is the analog of graphical deletion. Any collection of state loops that are connected by black nodes is evaluated as $\delta.$ Loops with common round nodes are regarded as disjoint so that $n$ loops with no black nodes will evaluate as $\delta^{n}.$
If we take $B=A^{-1}$ and $\delta = -A^2 - A^{-2},$ then the Chromatic Bracket is invariant under the second and third Reidemeister moves and, by dint of the square virtual evaluation formula
it is a Laurent polynomial valued invariant. Each picture has been translated into an invariant polynomial. We explain below in Figures~\ref{Coeffs},\ref{EFF9},\ref{EFF10} why the deletion-contraction expansion for square virtuals gives an evaluation that respects the detour moves\\

We will use the same notation $\langle K \rangle$ for both the general bracket and for the chromatic bracket evaluated at $B = A^{-1}$ and $\delta = -A^2 - A^{-2}$. The general bracket takes its values in sums of Laurent polynomials in $A$ and $A^{-1}$ as coefficients of graphical equivalence classes of purely virtual diagrams. There will be cases where the chromatic method is not sufficient and we will examine these equivalence classes directly or by using another way to evaluate them. This will be detailed below.\\

Figure~\ref{EFF6} illustrates our rule for the nodal crossing. For this purpose, the reader should view Figure~\ref{EFF7}. In that Figure we show the expansion rules for the chromatic bracket polynomial. Here we use coefficients $A$ and $B$ for the two smoothings and a value of $\delta$ for the loop value.\\

The chromatic bracket polynomial is defined for a theory with two virtual crossings that we will indicate with round and square indicators at the 
corresponding flat crossing. The polynomial has the usual bracket expansion at a crossing coupled with our expansion for the square virtual crossings as we have used them for the graph theory earlier in the paper. In the formulas below we use the topological coefficients for the bracket with 
$B = A^{-1}$ and $\delta = -A^2 - A^{-2}.$ The resulting evaluation is invariant under all moves except the first Reidemeister move where it changes by a multiple of $(-A^3)$ or $(-A^{-3})$ for positive and negative curls respectively.\\

The context for the chromatic  bracket is multi-virtuals with two crossing types: circle and square and the rule that the square virtual crossing is expanded into twice the nodal crossing minus the circle virtual. Hence one expands the bracket until there are only collecctions of virtual loops.
Then the square virtuals are expanded, and one has virtual loops with nodes. A collection of loops that are connected by nodes is evaluated as a single $\delta$. Thus if $S$ is a state in the expansion of this bracket then the evaluation of $S$, denoted $[S] = \delta^{c(S)}$ where $c(S)$ denotes the number of connected components of $S$ relative to the 
nodes. We illustrate this calculation on the virtual flat link in Figure~\ref{EFF8}.\\

\noindent {\bf Remark.} Motivation for the chromatic bracket in terms of graph theory can be found in \cite{MVKT}.\\

We have the basic theorem:\\

\noindent {\bf Theorem  \cite{MVKT}.} The generalized bracket is invariant under detour moves for any choice of commuting variables $x,y, \delta.$
The generalized and chromatic brackets $<K>$ are invariant under the Reidemeister moves $R2$ and $R3$ when we take $B= A^{-1}, \delta = -A^2 - A^{-2}.$
We call this the {\it topological specialization} of the bracket. Thus the topological specialization for the chromatic bracket is an invariant of 
double multi-virtual knots and links and is a well defined Laurent polynomial. For arbitrary multi-virtuals the topological specialization of the generalized bracket is invariant under the R2 and R3 moves and takes values in the ring generated by the flat virtual classes of the states of the diagram with coefficients in the Laurent polynomials in $A.$\\

\noindent {\bf Proof.} It is only necessary to check that the expansion of the square virtual is compatible with the detour moves. This is verified diagrammatically in Figure~\ref{EFF9} and Figure~\ref{EFF10}. This completes the proof. $\hfill\Box$


\begin{figure}
     \begin{center}
     \begin{tabular}{c}
     \includegraphics[width=10cm]{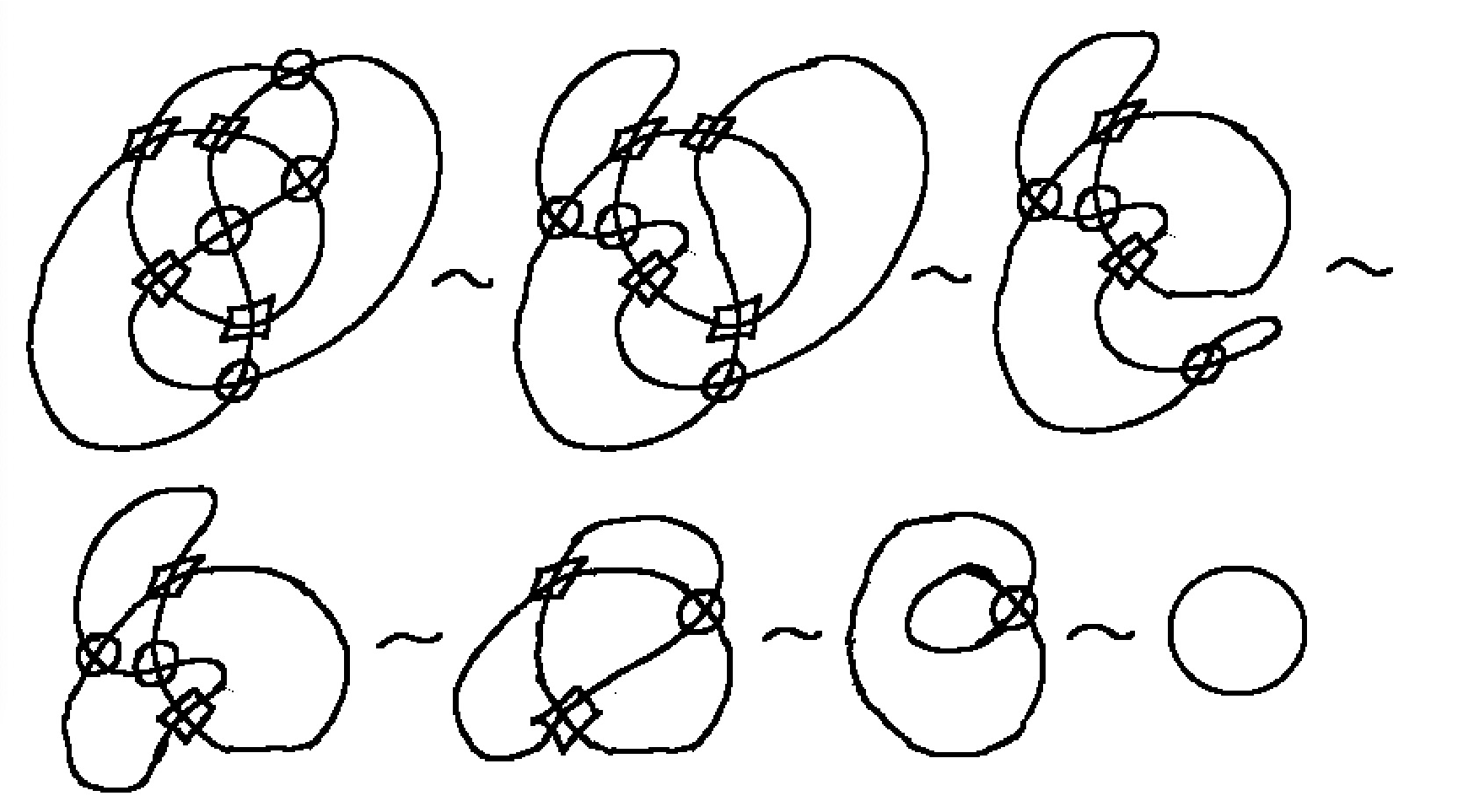}
     \end{tabular}
     \caption{\bf Single Loop with Two Virtual Crossings}
     \label{onecomp}
\end{center}
\end{figure}

\begin{figure}
     \begin{center}
     \begin{tabular}{c}
     \includegraphics[width=8cm]{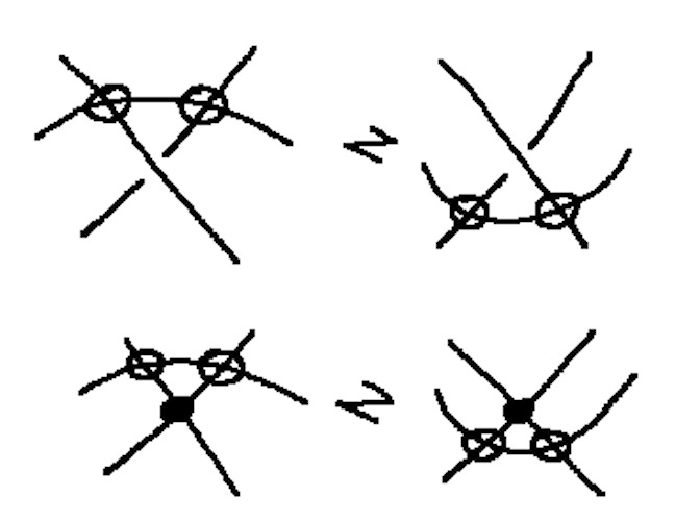}
     \end{tabular}
     \caption{\bf Detour Moves}
     \label{EFF4}
\end{center}
\end{figure}

\begin{figure}
     \begin{center}
     \begin{tabular}{c}
     \includegraphics[width=8cm]{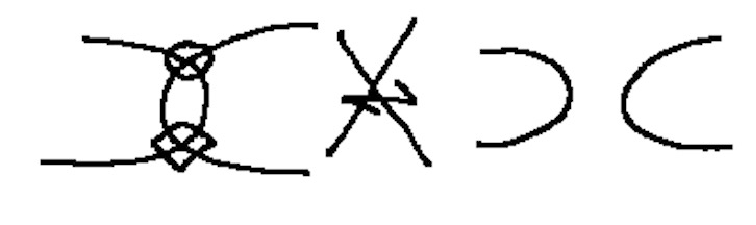}
     \end{tabular}
     \caption{\bf Non Cancellation of Distinct Virtual Crossings.}
     \label{EFF5}
\end{center}
\end{figure}

\begin{figure}
     \begin{center}
     \begin{tabular}{c}
     \includegraphics[width=8cm]{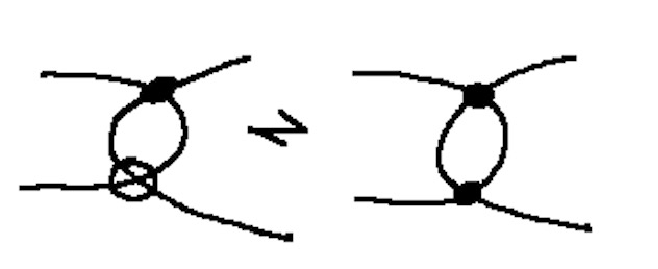}
     \end{tabular}
     \caption{\bf Nodal Equivalence}
     \label{EFF6}
\end{center}
\end{figure}

\clearpage

\begin{figure}
     \begin{center}
     \begin{tabular}{c}
     \includegraphics[width=8cm]{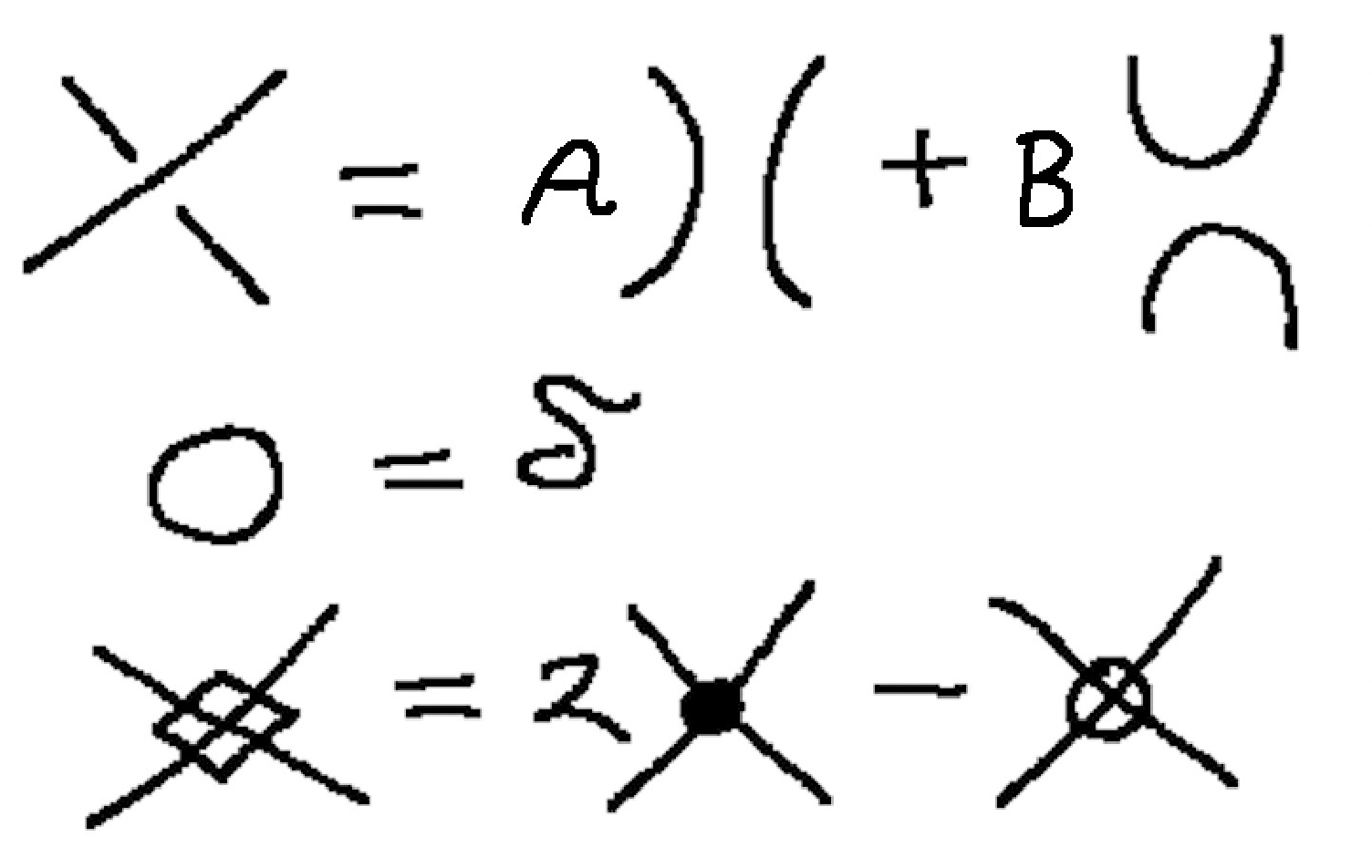}
     \end{tabular}
     \caption{\bf Chromatic Bracket}
     \label{EFF7}
\end{center}
\end{figure}

In Figure~\ref{EFF8} we illustrate double virtual flat calculations for the simple link using two virtual crossings, and another example with three components.\\

\begin{figure}
     \begin{center}
     \begin{tabular}{c}
     \includegraphics[width=8cm]{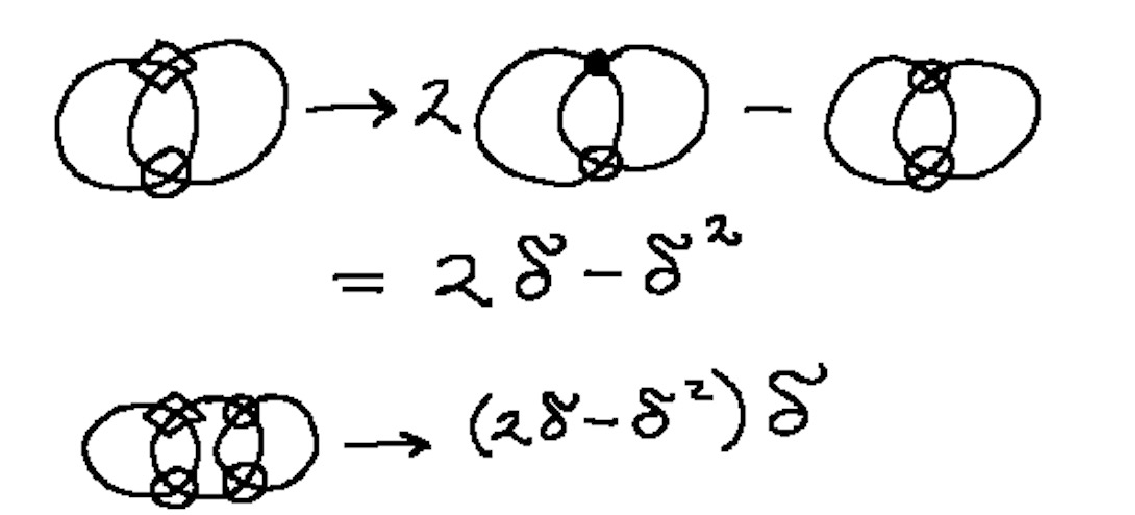}
     \end{tabular}
     \caption{\bf Loop Evaluations}
     \label{EFF8}
\end{center}
\end{figure}

In Figure~\ref{Coeffs} shows that the expansion for the square virtual crossing necessarily has the coefficients $2$ and $-1$ if it is to satisfy the simplest flat Reidemeister two detour move.
In Figure~\ref{EFF9} and Figure~\ref{EFF10} we illustrate the graphical calculations that verify that the expansion formula for the square virtual crossing is compatible with 
the detour moves and flat Reidemeister three moves (also detour moves) for the doubled virtual crossings.  These figures constitute the crucial lemma of  \cite{MVKT} that this method of evaluation leads to a regular isotopy invariant for doubled virtual links.\\

\begin{figure}
     \begin{center}
     \begin{tabular}{c}
     \includegraphics[width=10cm]{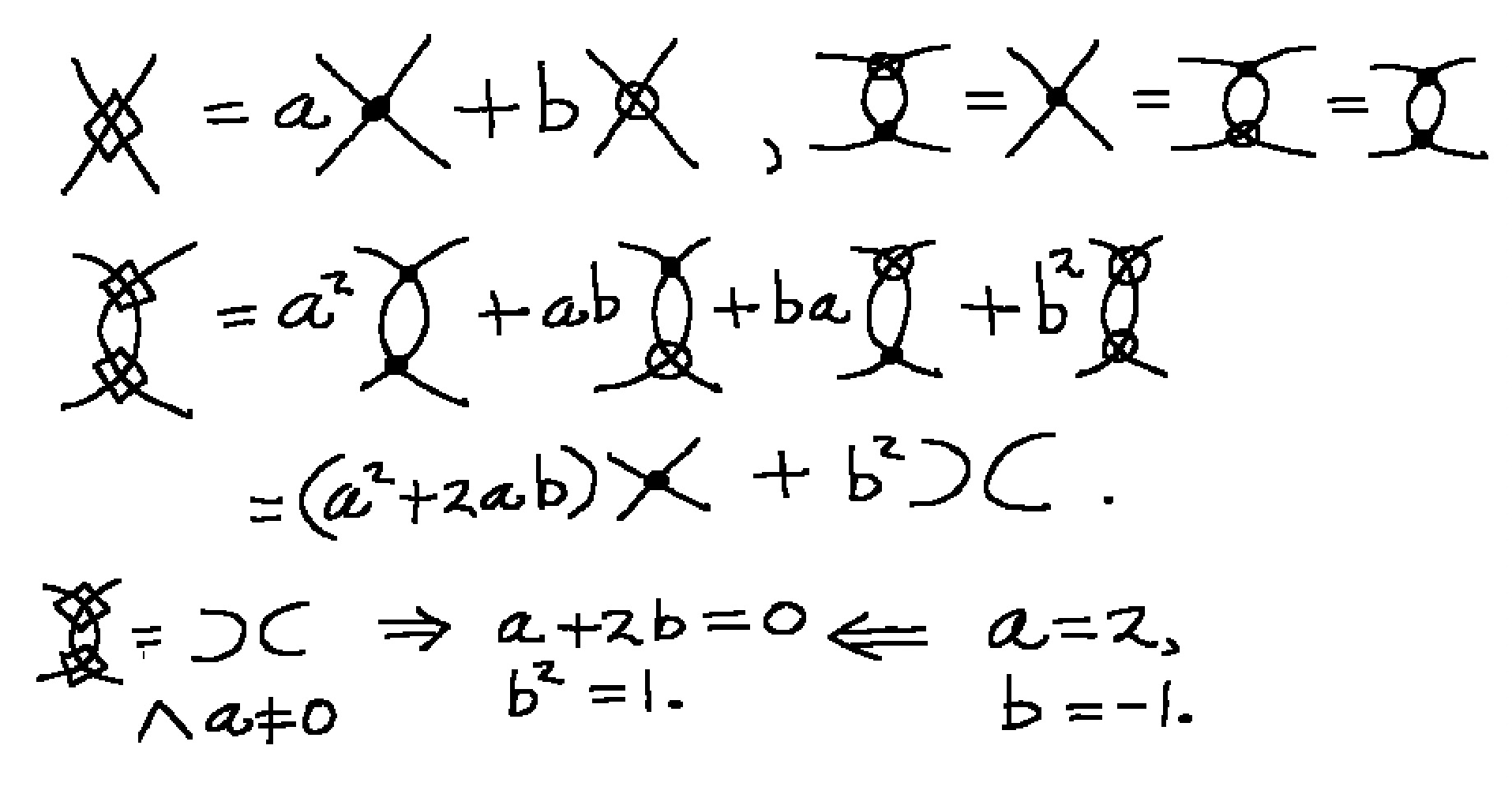}
     \end{tabular}
     \caption{\bf Graphical Expansion Verification}
     \label{Coeffs}
\end{center}
\end{figure}

\begin{figure}
     \begin{center}
     \begin{tabular}{c}
     \includegraphics[width=10cm]{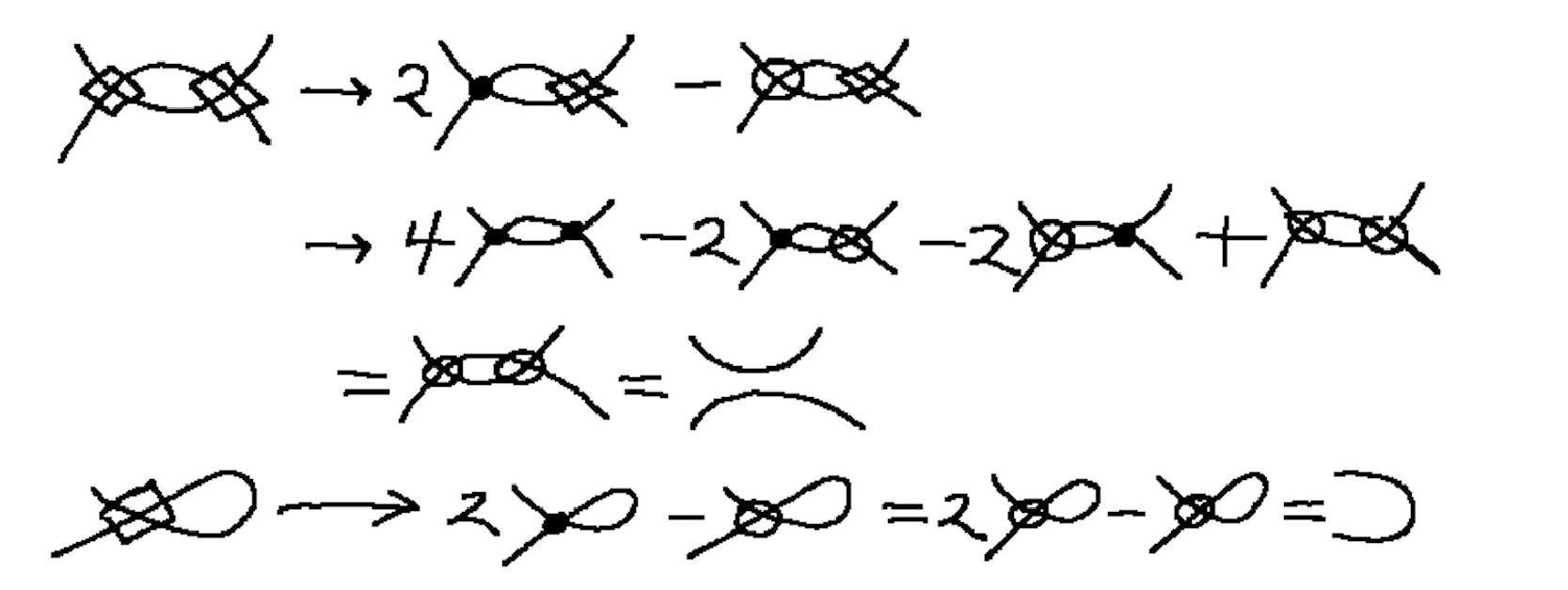}
     \end{tabular}
     \caption{\bf Graphical Expansion Verification}
     \label{EFF9}
\end{center}
\end{figure}

\begin{figure}
     \begin{center}
     \begin{tabular}{c}
     \includegraphics[width=10cm]{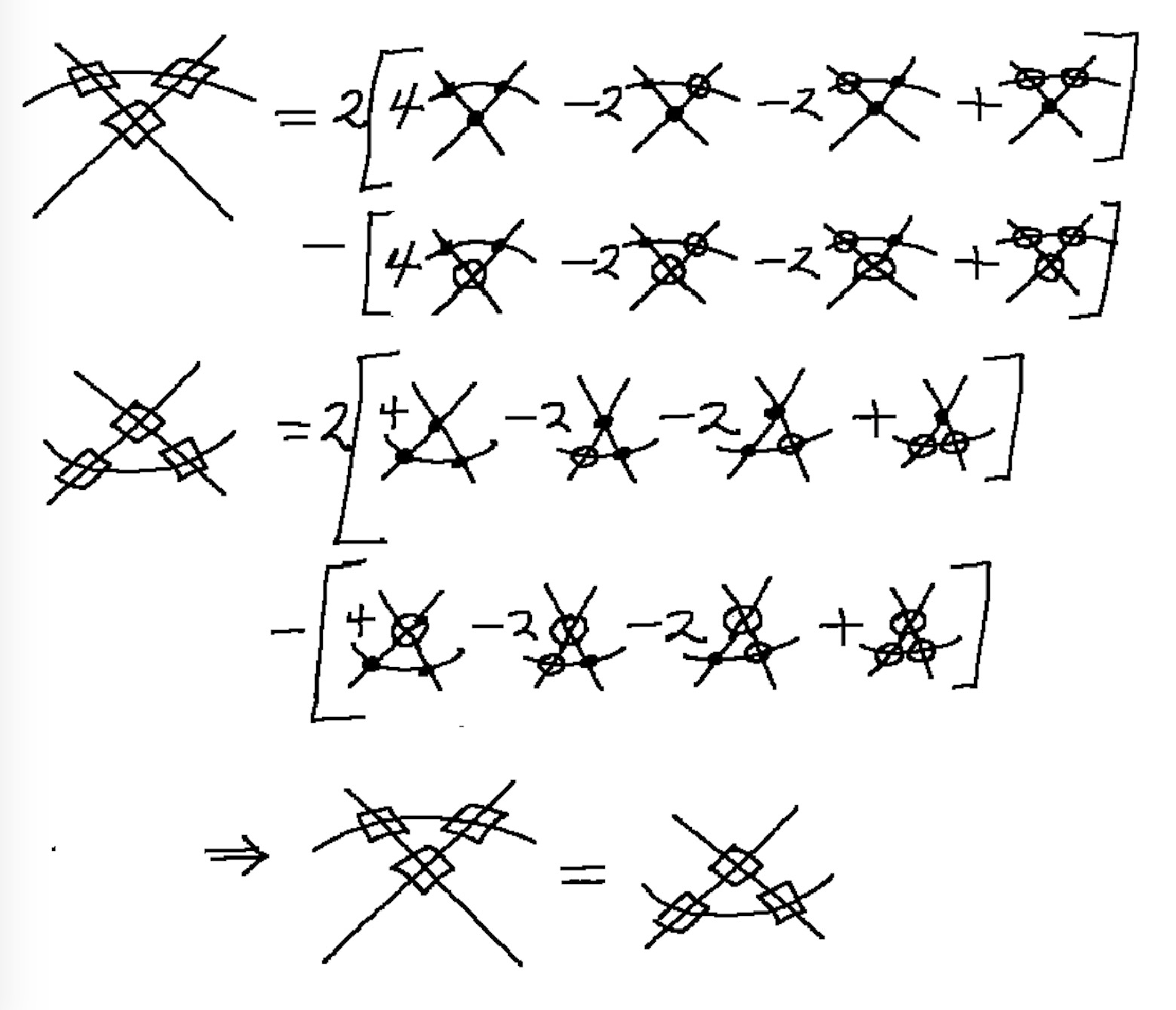}
     \end{tabular}
     \caption{\bf Graphical Expansion Verification}
     \label{EFF10}
\end{center}
\end{figure}

\clearpage

\begin{figure}
     \begin{center}
     \begin{tabular}{c}
     \includegraphics[width=10cm]{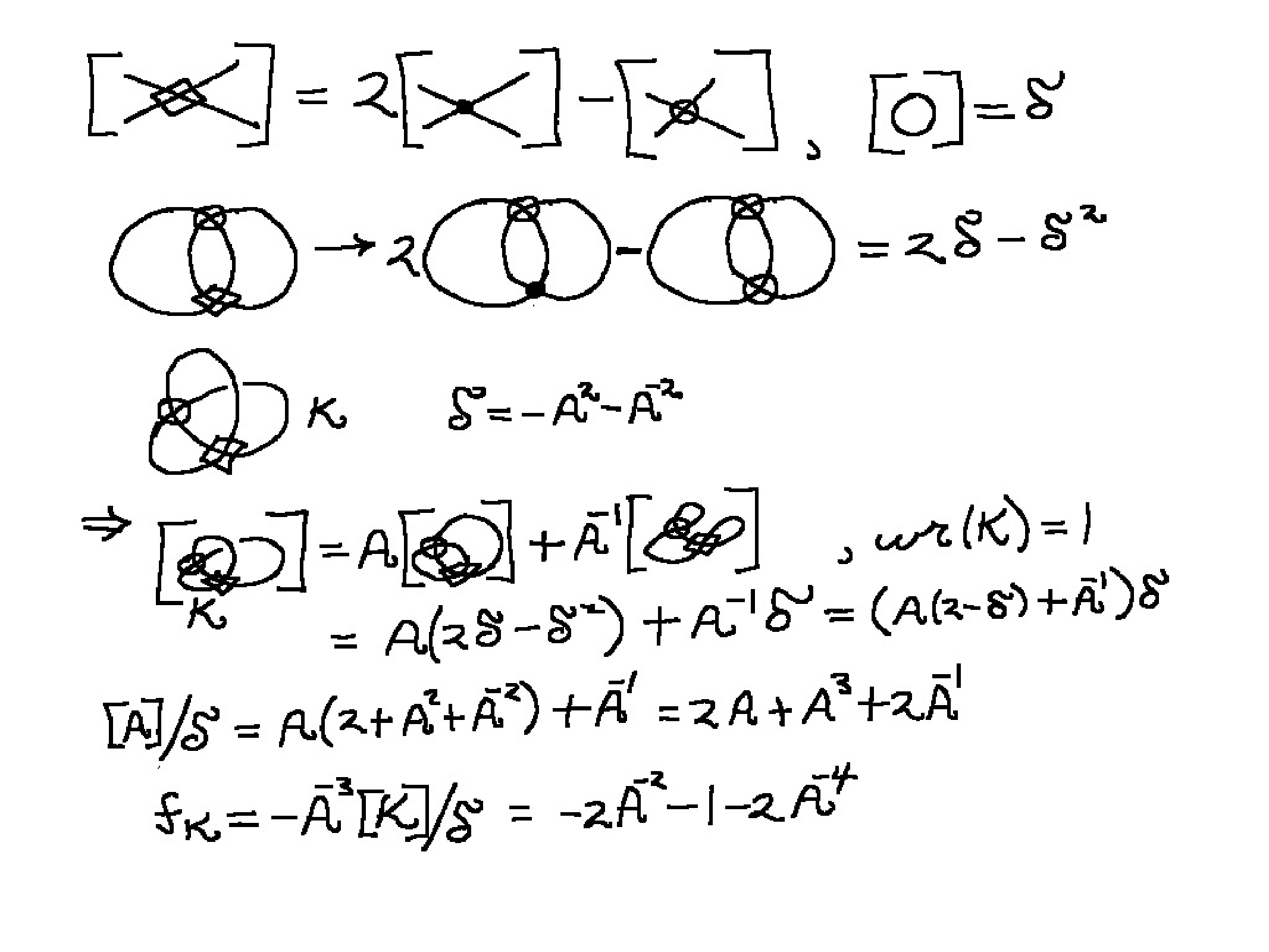}
     \end{tabular}
     \caption{\bf A Chromatic Bracket Evaluation}
     \label{EFF12}
\end{center}
\end{figure}

Figure~\ref{EFF12} shows the expansion of the generalized bracket for double virtuals on a trefoil diagram with one classical crossing, one square virtual crossing and one round virtual crossing.\\

\begin{figure}
     \begin{center}
     \begin{tabular}{c}
     \includegraphics[width=12cm]{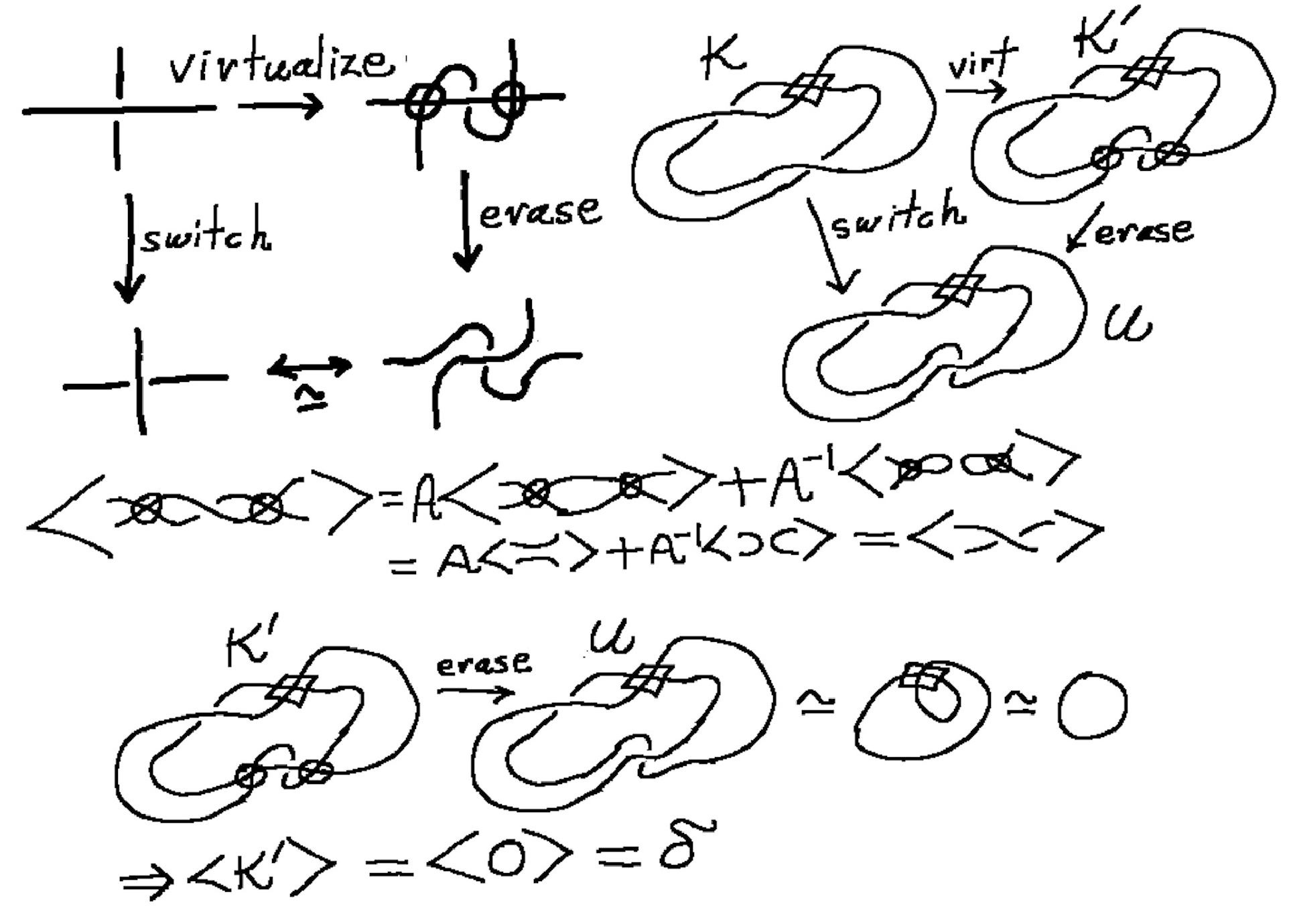}
     \end{tabular}
     \caption{\bf Virtualization for Multi-Virtuals}
     \label{virt}
\end{center}
\end{figure}

\begin{figure}
     \begin{center}
     \begin{tabular}{c}
     \includegraphics[width=12cm]{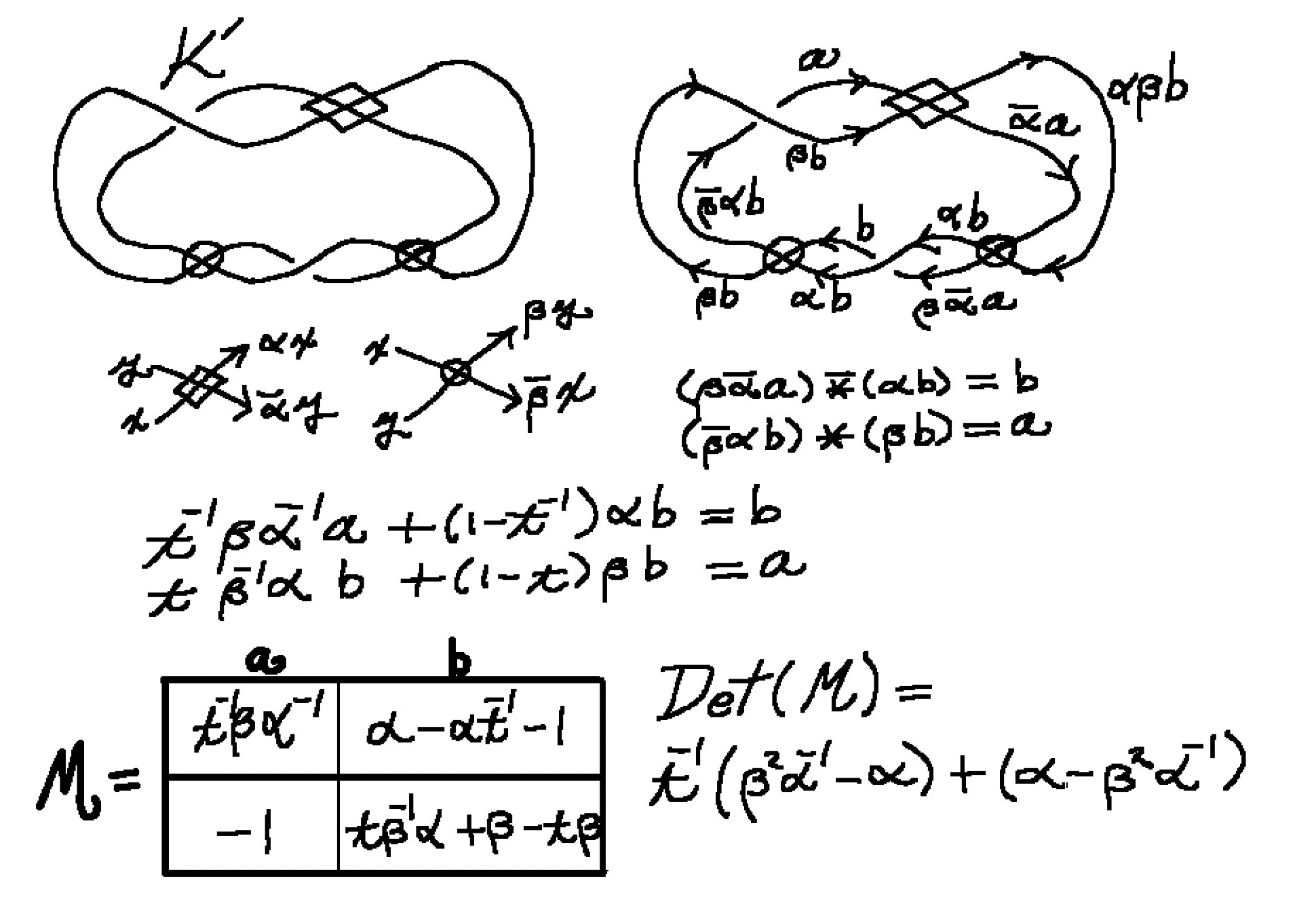}
     \end{tabular}
     \caption{\bf Non-Triviality of a Virtualized Knot}
     \label{nontriv}
\end{center}
\end{figure}

From Figure~\ref{virt} onward in this paper, we begin new material that has not been discussed in pervious work for multi-virtuals.\\

Figure~\ref{virt} shows how the {\it virtualization operation} changes a crossing by reversing its undercrossing direction and flanking it by two virtuals of the same type.
This operation can be described by taking a walk on the diagram and reversing the direction in which one goes under the crossing by first going virtually across the overcrossing arc, then under the crossing, then back across the overcrossing 
arc. The upper left of the figure shows that {\it erasure of the virtual crossings in the virtualization results in the switch of the original crossing}. Just below that part of the figure, we show that the generalized bracket evaluation of 
 a diagram with a virtualized crossing is identical with the evaluation of the bracket with the erasure of the virtuals flanking the crossing. Thus the bracket evaluation of a diagram with a virtualized crossing is the same as the evaluation
 of the diagram obtained by switching that crossing. It can happen that switching a crossing can produce the diagram of an unknot. in such a case, the virtualized diagram will have a trivial bracket polynomial. We illustrate this situation in the figure with the original knot $K$ and its virtualization $K'.$ Thus the generalized bracket of $K'$ is the same as the bracket of an unknot. $K'$ is in fact a non-trivial and non-classical multivirtual knot. To verify this, see Figure~\ref{nontriv}. There we have used
 the generalized quandle technique of \cite{MVKT} to compute a multi-virtual generalized Alexander polynomial for $K'.$ The calculation shows that $K'$ is non-trivial and non-classical. We refer the reader to the paper just mentioned for the 
 details of this technique.\\
 
 \noindent {\bf Remark.} The virtualization method produces infinitely many non-trivial non-classical virtual and multi-virtual knots with trivial bracket polynomial (unit Jones polynomiall). In the case of single virtual crossing theory, this situation has been extensively investigated \cite{vkt,DVK,DKK}. We know that virtual knots obtained from non-trivial classical knots by virtual moves and the virtualization process applied to some crossings (so that the switched classical diagrams are unknots) are non-classical and non-trivial virtual knots. More information may emerge from the multi-virtual instances of knots with trivial Jones polynomial.\\

\section{Strict Virtual Linkoids}

We now discuss the definition of multi-virtual linkoids in the strict sense and how closures provide invariants by using either classical crossings or virtual crossings. The key idea for new invariants is to use new virtual crossings outside the set of virtual crossings that are used for the linkoids.\\ 

\begin{figure}
     \begin{center}
     \begin{tabular}{c}
     \includegraphics[width=12cm]{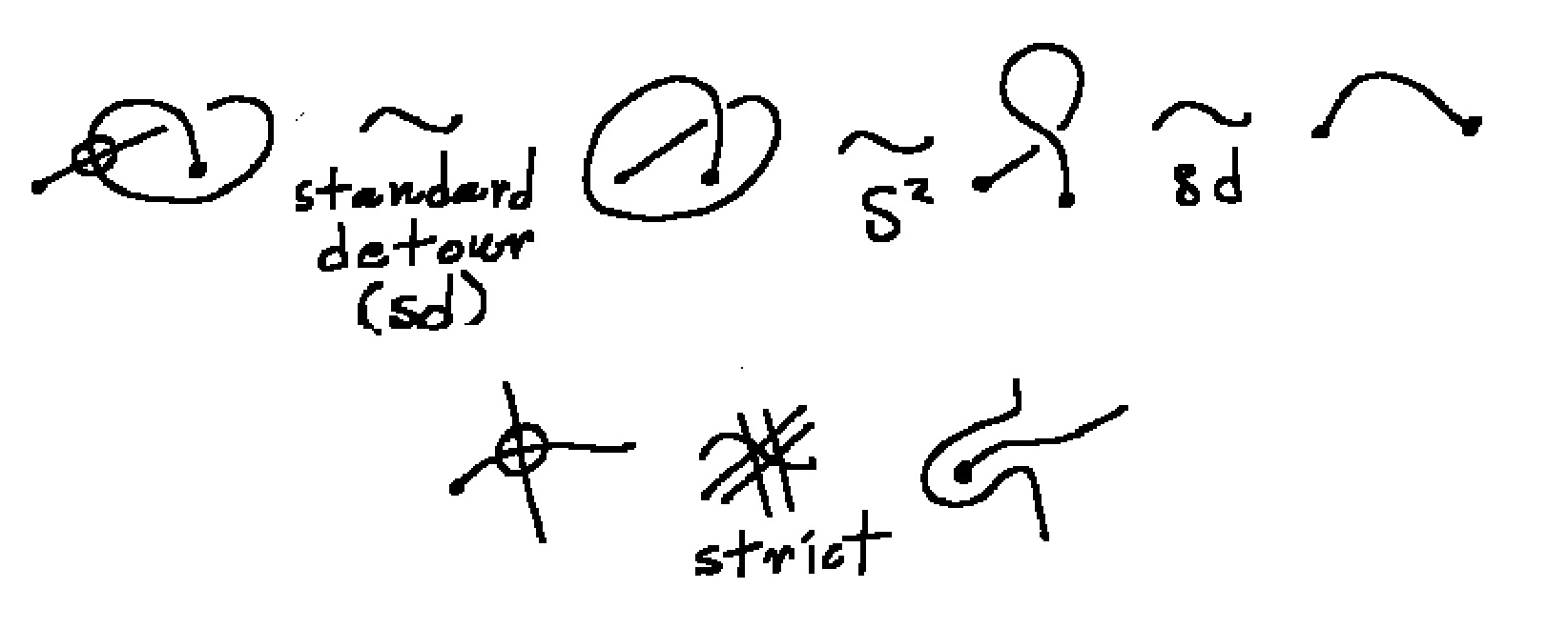}
     \end{tabular}
     \caption{\bf Strict and Standard Equivalence}
     \label{strict}
\end{center}
\end{figure}

A {\it knotoid} is a tangle diagram with two ends so that the endpoints are in possibly distinct regions \cite{Turaev,Neslihan}. See Figure~\ref{strict} for an illustration of a virtual knotoid diagram. In the original papers on knotoids
one is concerned with diagrams for knotoids that have only classical crossings. In this standard theory one allows classical Reidemeister moves and one forbids any move that would pass an arc in the diagram across an endpoint.
Thus the standard knotoids are very similar to classical knot diagrams except for the fact that the endpoints of the diagram can be in different regions in the plane. We can consider {\it linkoids} as well where there are many components. Usually in referring to a linkoid one has still in mind that there are two endpoints in the diagram. If we extend further to diagrams with more than two endpoints, these are called {\it multi-linkoids}. It is natural to have virtual knotoids, linkoids and multi-linkoids. And we can have multi-virtual species of all of these with diagrams containing more than one virtual crossing type. Diagrams will be taken to be on the surface of the two-sphere $S^2.$ The key point that we now address is the question of dealing virtual moves for knotoids and their generalizations.\\

\begin{figure}
     \begin{center}
     \begin{tabular}{c}
     \includegraphics[width=12cm]{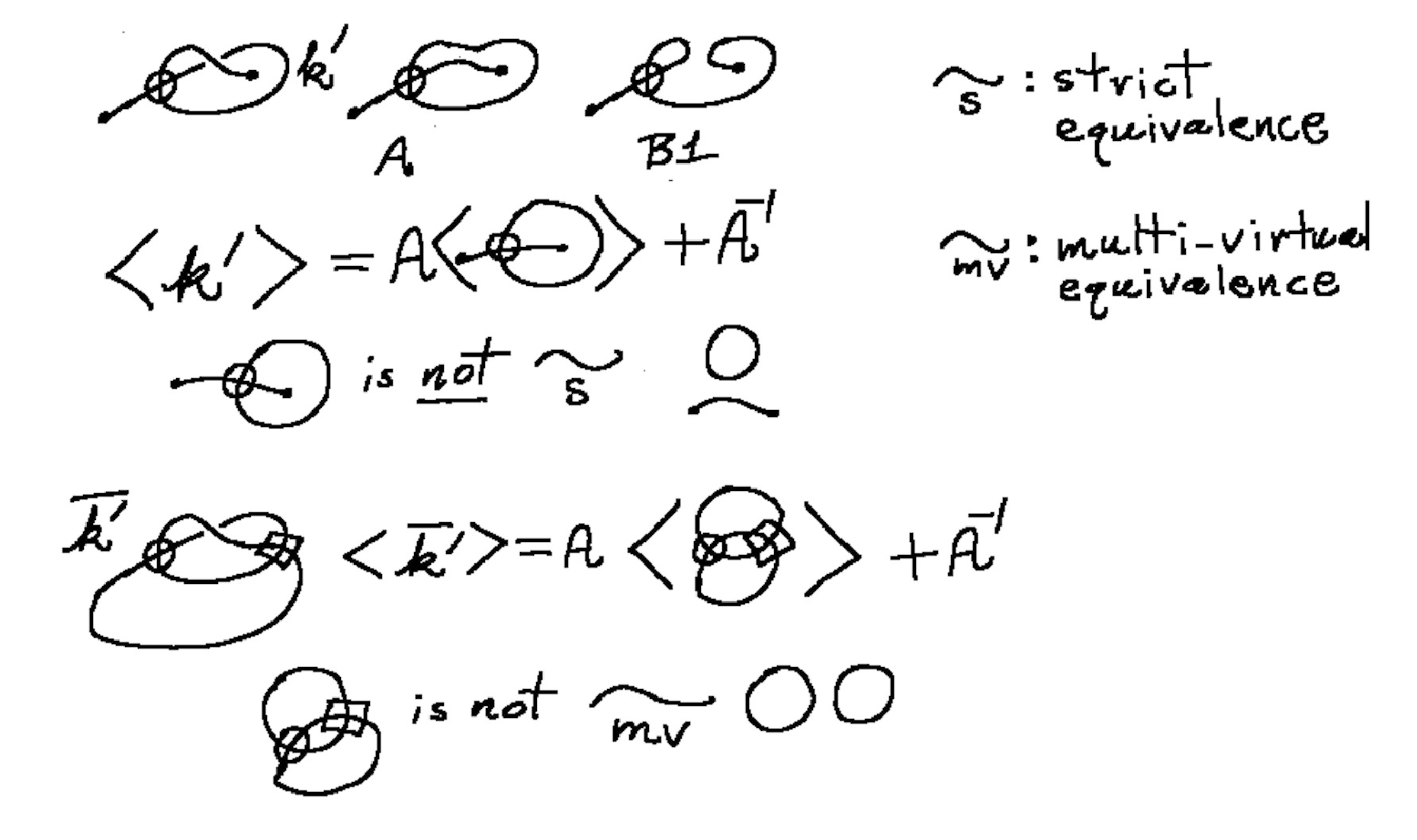}
     \end{tabular}
     \caption{\bf Strict Equivalence and Multi-Virtual Closure}
     \label{closure}
\end{center}
\end{figure}

\begin{figure}
     \begin{center}
     \begin{tabular}{c}
     \includegraphics[width=12cm]{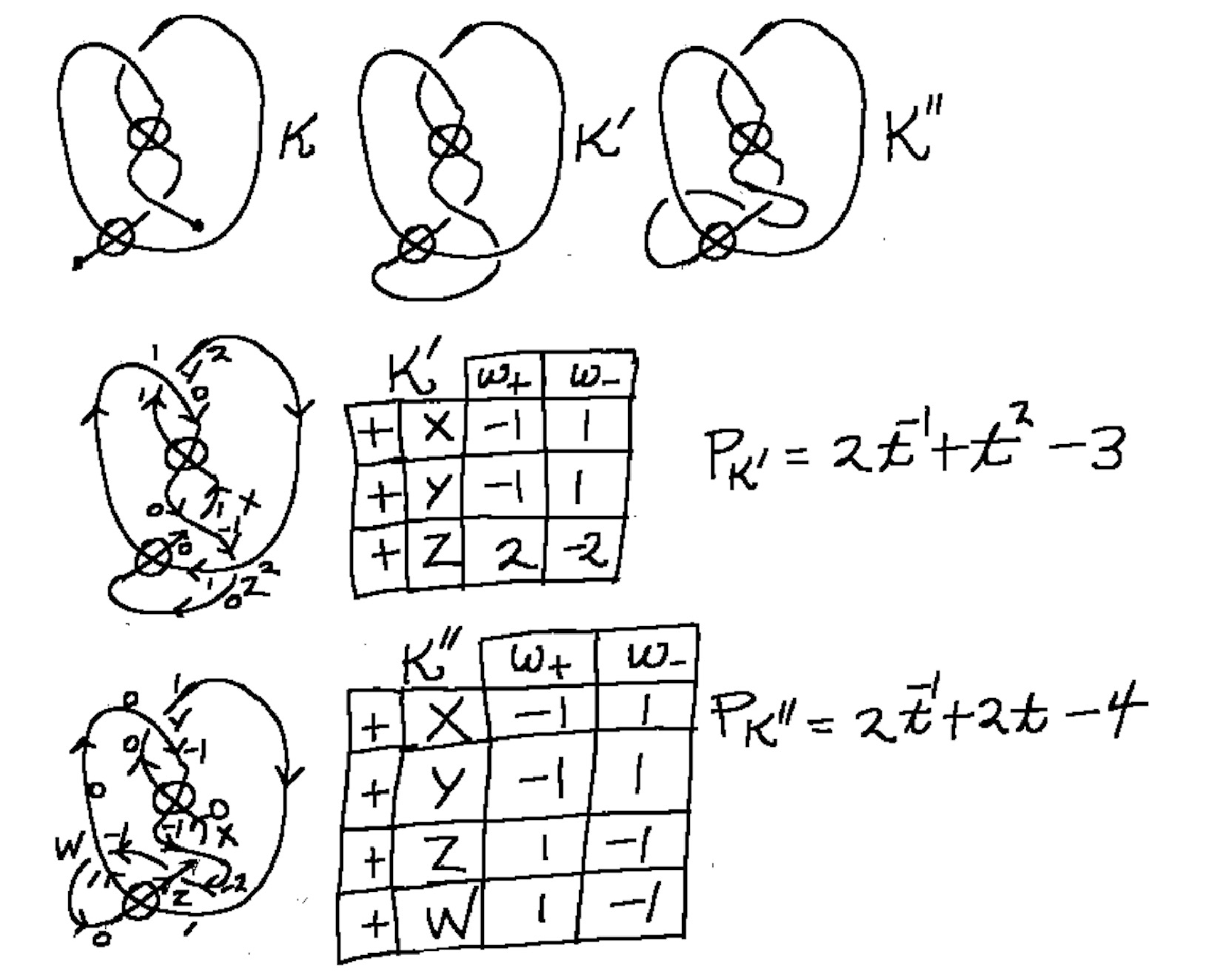}
     \end{tabular}
     \caption{\bf Classical Underclosure is Ill-Defined for Virtual Linkoids}
     \label{underclosure}
\end{center}
\end{figure}

\begin{figure}
     \begin{center}
     \begin{tabular}{c}
     \includegraphics[width=12cm]{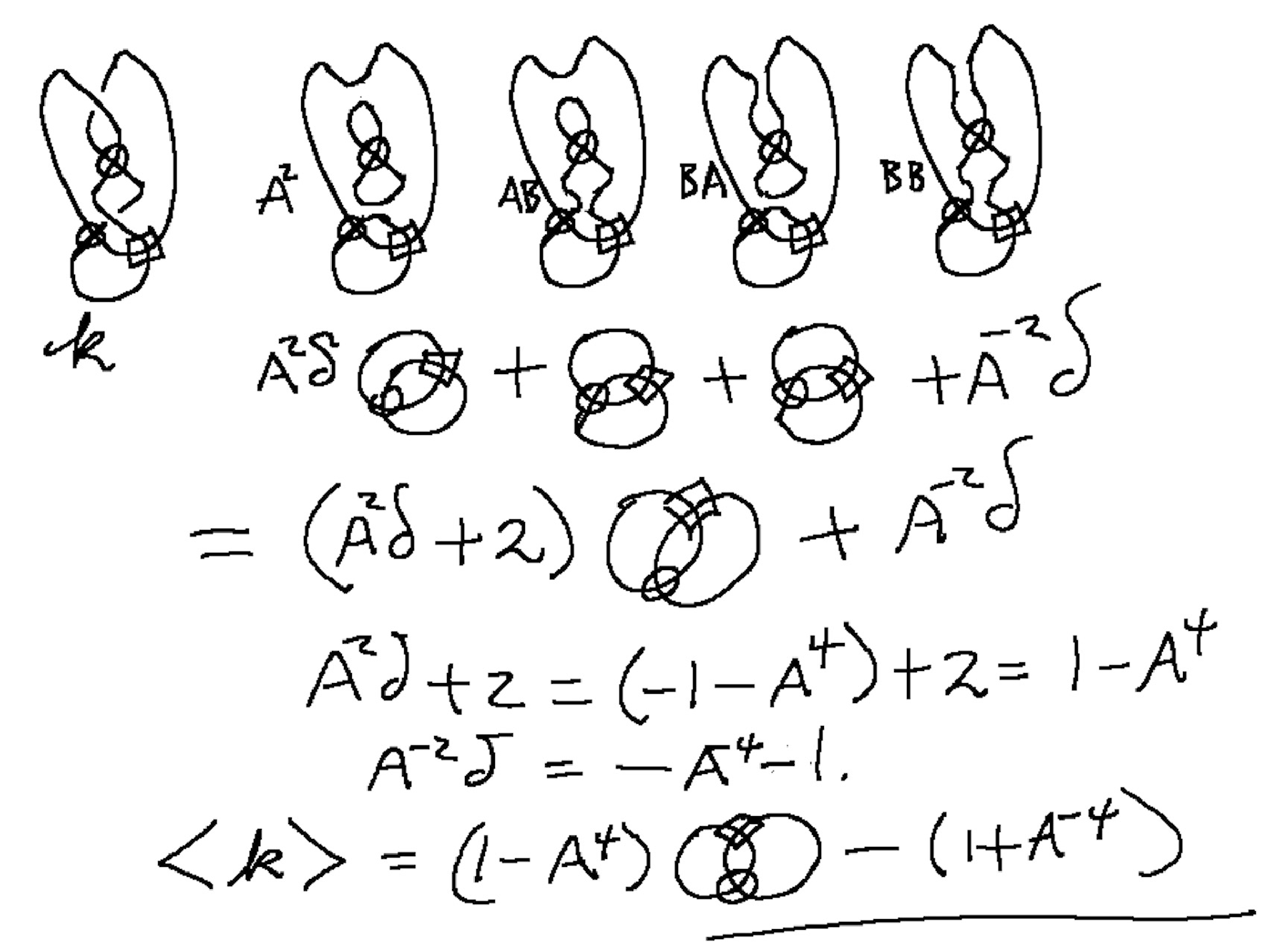}
     \end{tabular}
     \caption{\bf A Multi-Bracket Evaluation}
     \label{multi}
\end{center}
\end{figure}

\begin{figure}
     \begin{center}
     \begin{tabular}{c}
     \includegraphics[width=12cm]{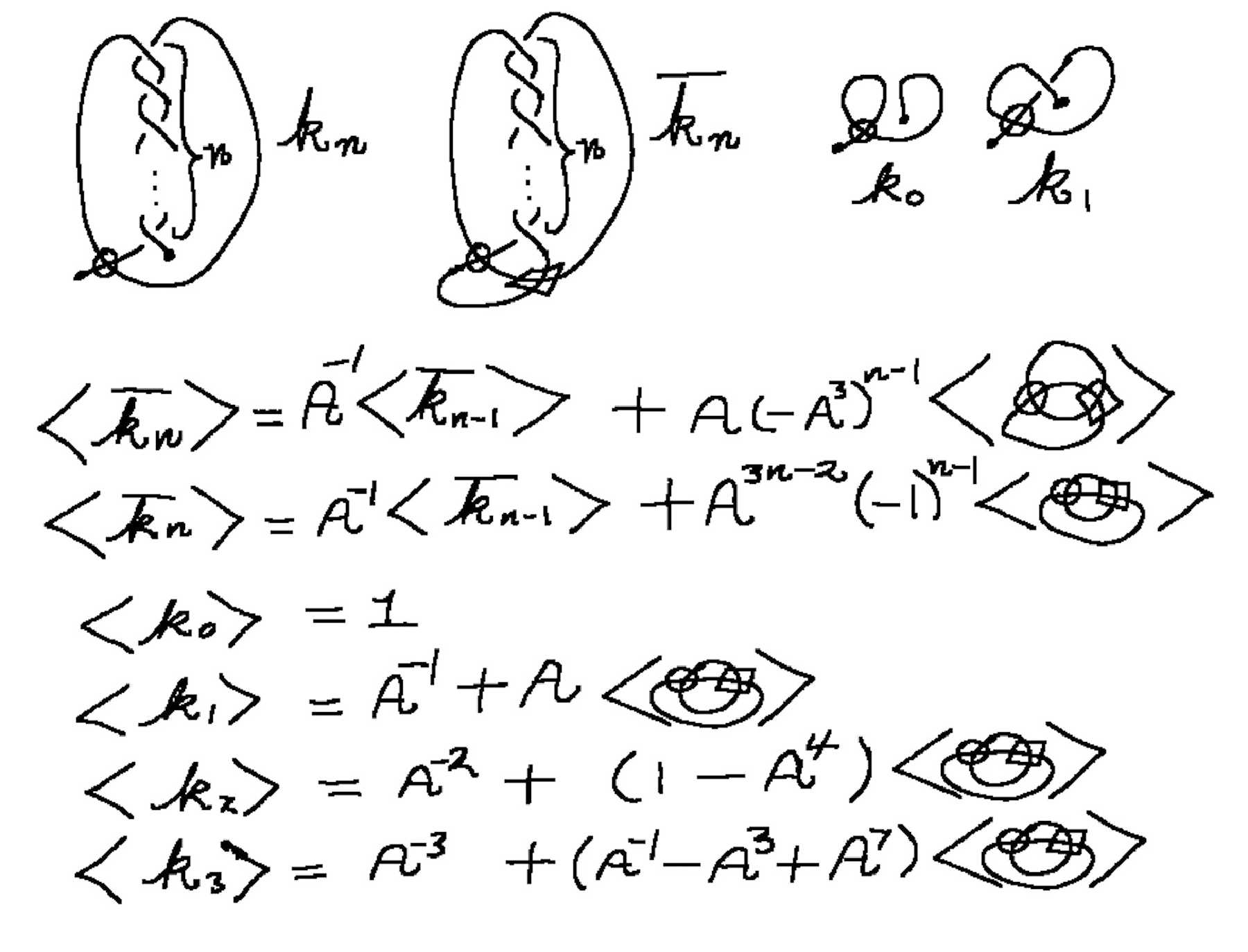}
     \end{tabular}
     \caption{\bf A family of linkoids and their closures.}
     \label{manymulti}
\end{center}
\end{figure}

\begin{figure}
     \begin{center}
     \begin{tabular}{c}
     \includegraphics[width=12cm]{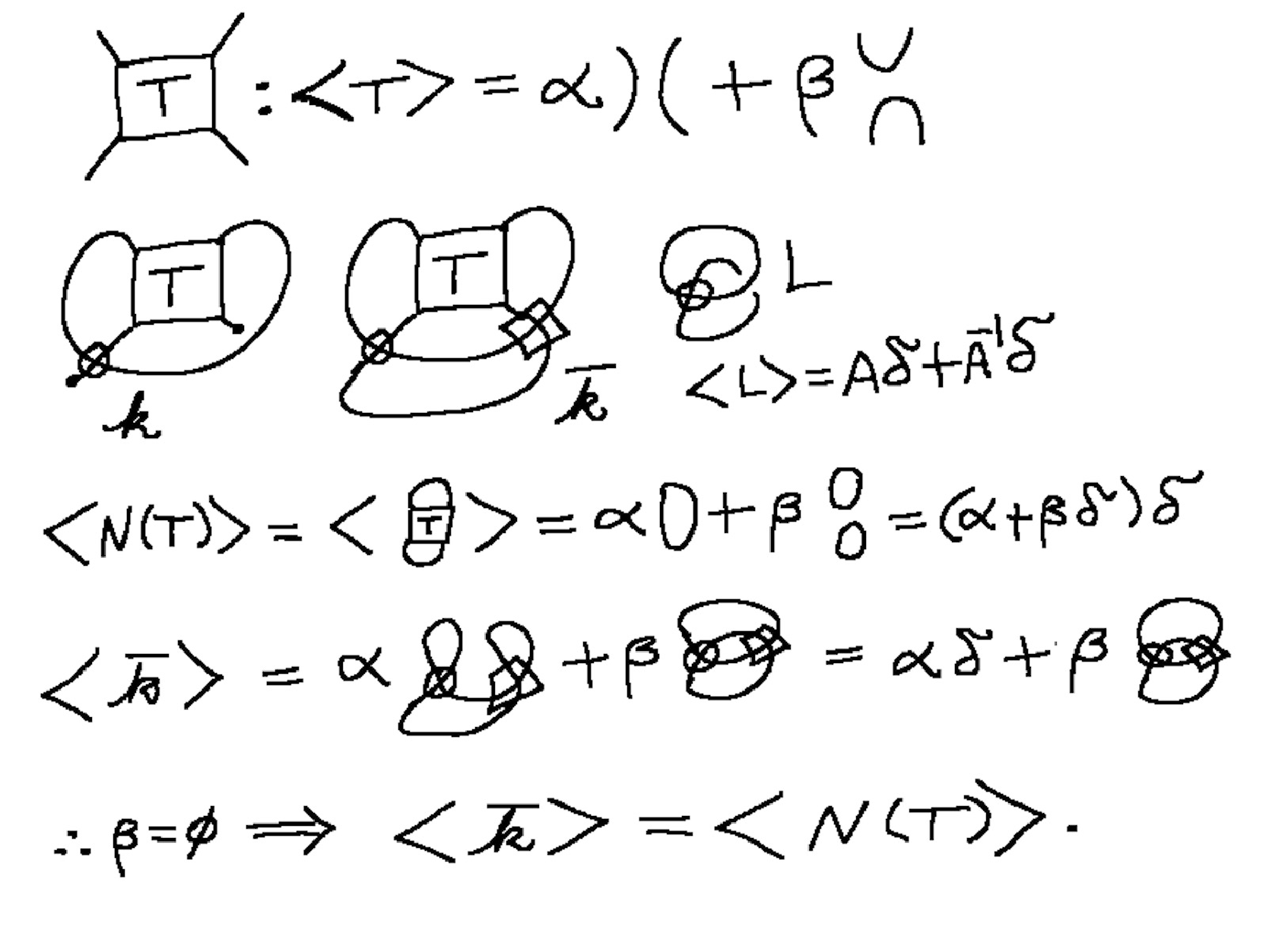}
     \end{tabular}
     \caption{\bf A family of linkoids defined by tangles.}
     \label{tangles}
\end{center}
\end{figure}

\begin{figure}
     \begin{center}
     \begin{tabular}{c}
     \includegraphics[width=12cm]{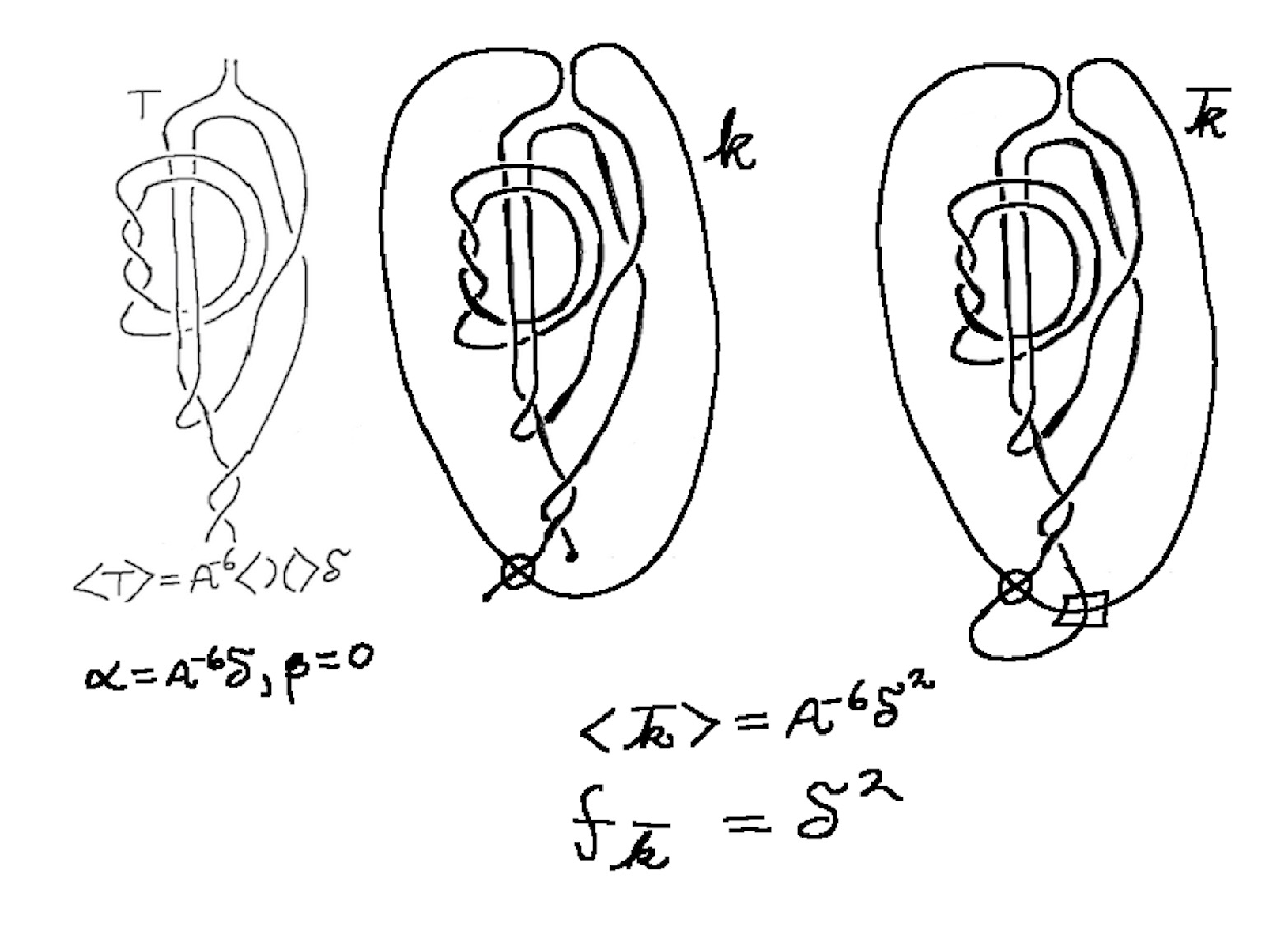}
     \end{tabular}
     \caption{\bf Producing examples not detectable by the bracket.}
     \label{invisible}
\end{center}
\end{figure}

\begin{figure}
     \begin{center}
     \begin{tabular}{c}
     \includegraphics[width=12cm]{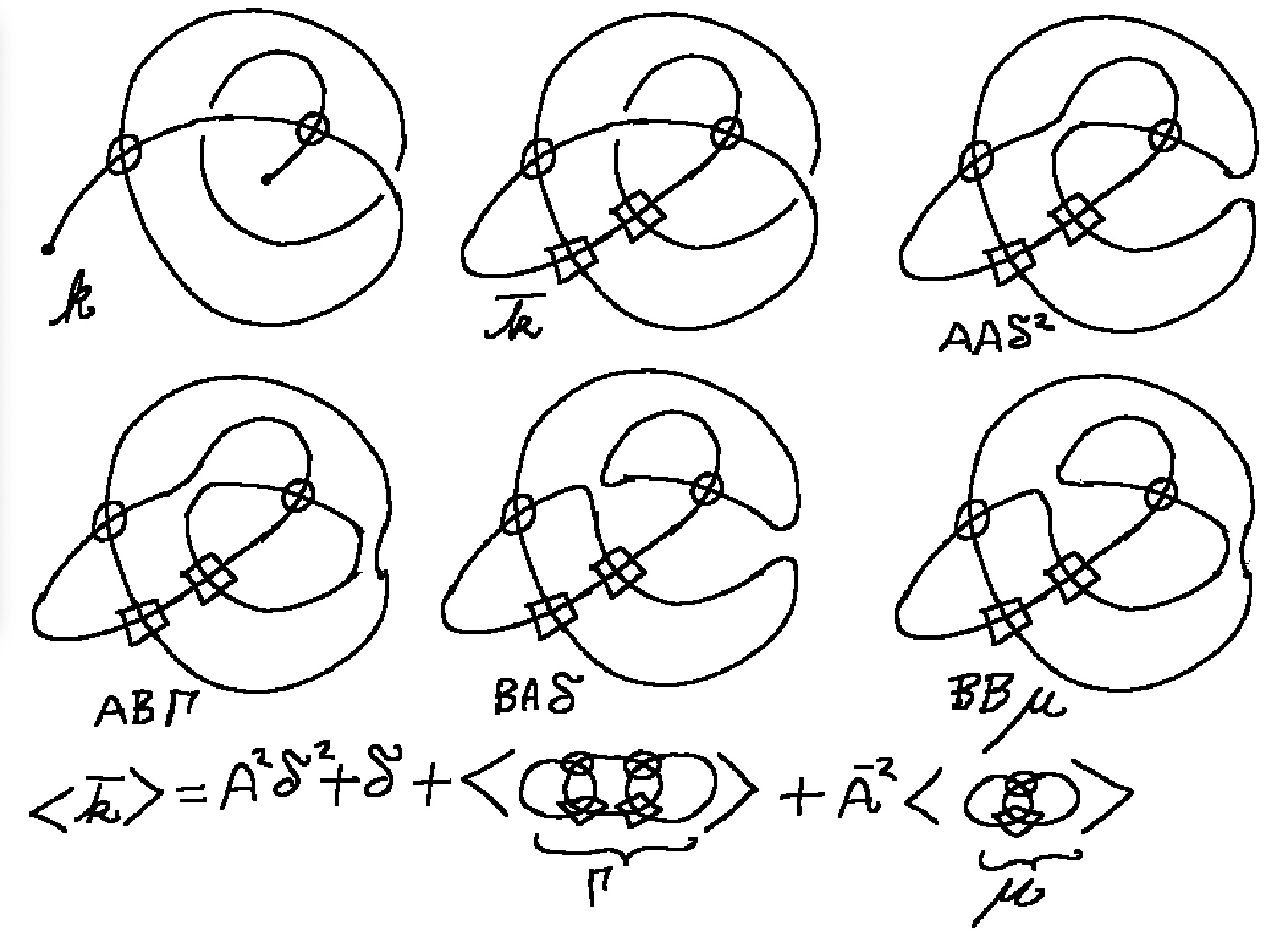}
     \end{tabular}
     \caption{\bf Using closure to estimate the height of a strict linkoid.}
     \label{height}
\end{center}
\end{figure}

Figure~\ref{strict} illustrates that unrestricted detour moves can unknot an example knotoid diagram. We define {\it strict} equivalence by using local Reidemeister type virtual moves as shown in 
Figure~\ref{Figure 1} and Figure~\ref{EFF3}. Restriction to such moves means that we can enforce the rule that {\it no move shall transfer an arc across an endpoint}.  
This defines strict equivalence of virtual and multi-virtual knotoids,linkoids and multi-linkoids.\\

In Figure~\ref{closure} we illustrate a generalized bracket expansion for a virtual knotoid $k'$ and a generalized bracket expansion for the ``box" virtual closure $\bar{k'}$ of $k'.$
The box virtual closure is obtained by connecting the endpoints of the knotoid $k'$ by an arc with square-shaped virtual crossings that differ from the ``round" virtual crossings that
are allowed for $k'.$ In general, the virtual closure of a knotoid will make a closure arc using a virtual crossing type that is not already present in the given knotoid. In the first bracket expansion we see that there is a term
consisting of a loop intersected by an arc with a single virtual crossing. This is a state in the state expansion of this generalized bracket, and it is a non-trivial state in regard to strict equivalence of states -- generated by strict detour moves.
We leave it as an exercise for the reader to prove that that this state is inequivalent to a state with an arc and a disjoint loop. Thus the first bracket expansion serves to prove the strict non-triviality of $k'.$ Nevertheless, we have a proof of strict non-triviality from the bracket expansion of the box closure $\bar{k'}$ since there we see the already analysed pair of loops intersecting in round and square virtual crossings. Thus this figure illustrates two techniques that can be used to detect strict non-triviality.\\

From now on in the paper when we refer to virtual knotoids, this includes all the variations: linkoids and multi-knotoids, unless otherwise specfied.\\

In Figure~\ref{underclosure} we illustrate how classical underclosure (and by the same token, overclosure) is ill-defined for virtual knotoids. The underclosure and overclosure operations are well-defined for knotoids 
with only classical crossings and are accomplished by connecting corresponding endpoints by arcs that are either beneath the plane (for the underclosure) or above the plane (for the overclosure). In extending such closures to virtual knotoids one finds that the type of the resulting virtual knot depends on the path of the arc. This dependence is due to the fact that virtual knots are not invariant under the forbidden moves of Figure~\ref{Figure 3}. In Figure~\ref{underclosure} we illustrate the computation of the affine index polynomial \cite{affineindex} for two choices of underclosrure arc in a given example of a knotoid. That the polynomials are distinct proves that the virtual knots obtained by underclosure are not equivalent. The reader can find the details of the definition of the affine index polynomial in the paper \cite{affineindex}.\\

In Figure~\ref{multi} we illustrate the box virtual closure of the virtual knotoid shown in Figure~\ref{underclosure}. Taking the virtual closure as a multi-virtual knot this knot is an invariant of the strict virtual class of the knotoid. 
The figure shows the generalized bracket calculation of this multi-virtual closure and shows via the calculation that the mult-virtual knot closure is non-trivial since it contains an irreducible two-virtual state.\\

Figure~\ref{manymulti} illustrates a family of virtual knotoids $k_n$ (strictly non-trivial for $n=1,2,3,\cdots,...$) and indicates the proof that they are are strictly non-trivial and strictly distinct via calculations of the bracket polynomials of their  box virtual 
closures.\\

Figure~\ref{tangles} indicates a family of virtual knotoids defined by tangles $T.$ Assuming these are classical tangles and that the expansion of the bracket of the tangle $T$ has the form 
$\langle T \rangle = \alpha \langle \infty \rangle + \beta \langle 0 \rangle,$ (where $\infty$ and $0$ are the ``infinity" and ``zero" tangles shown in the figure). the figure shows how to calculate the box virtual closures of these knotoids and concludes that their strict non-triviality is detected when $\beta \ne 0.$\\

In Figure~\ref{invisible} we give an example of a linkoid that whose strict class is not seen to be different from a disjoint union of an arc and a circle by the generalized bracket. The method employed here is to use a tangle as in the previous example such that
$\beta = 0.$ This tangle produces a classical link (by standard numerator closure $N(T)$ of the tangle as in Figure~\ref{tangles}) that is invisible under the bracket from an unliink of two components. See \cite{EKT}, We have modified this construction to give linkoids that are invisible to the generalized bracket. In a subsequent paper we will give a proof that this linkoid is strictly non-trivial.\\

In Figure~\ref{height} we illustrate how to use the box virtual closure of a knotoid to determine its {\it height}. The height of a knotoid is the least ``distance" between its endpoints among all diagrams of that knotoid.
This distance is the number of arcs of the knotoid that need be crossed in order to connect the endpoiints. We have the following lemma.\\

\noindent {\bf Lemma.} Let $k$ be strict knotoid with two endpoints. Let $\bar{k}$ denote the box virtual closure of $k$ so that $\bar{k}$ is a multivirtual knot with round and box virtual crossings. Then among all diagrams of $\bar{k}$ , up to multivirtual equivalence, there are those with the least number of box virtual crossings. Let $S(k)$ denote the least number of box virtual crossings among all diagrams of $\bar{k}.$ Then $S(k)$ is a lower bound for the height of $k$ as a strict knotoid.
Furthermore, the maximal number of appearances of the box virtual crossing in the reduced states of the bracket $\langle \bar{k} \rangle$ is a lower bound for $S(k).$ Thus one can use the generalized bracket to estimate the height of a strict knotoid.\\

\noindent {\bf Proof.} The proof follows at once from the definitions of the objects and invariants discussed in the statement of the Lemma. $\hfill\Box$\\

Using the Lemma and the calculation shown in Figure~\ref{height}, we see that the strict knotoid depicted in that figure has height equal to 2.\\

This section has been a selection of examples showing how the box virtual closure of a knotoid can be used to determine information about its strict knotoid classification. \\

\section{Planar Strict Linkoids and Generalized Loop Bracket Polynomials}
In this section we consider strict linkoids in the plane. We are already familiar \cite{Turaev,Neslihan} with the differences between linkoids on the two-sphere and linkoids in the plane.
When working with the bracket polynomial for linkoids in the plane it is possble for a loop in the state explansion to encircle an arc with endpoints of the linkoid. Under such circumstances the loop will remain encircling the arc throughout all
isotopies of the linkoid and its state expansion. Thus Turaev \cite{Turaev} defined the {\it loop bracket polynomial} for knotoids to include variables$U_{n}$ that correspond to $n$ nested circles surrounding a knotoid arc. Similar extensions can be done for linkoids as in \cite{Neslihan}. Here we consider multi-virtual linkoids in the plane under strict equivalence and show phenomena that occur in their states. These are phenomena that generalize the special variables of the loop bracket. \\

\begin{figure}
     \begin{center}
     \begin{tabular}{c}
     \includegraphics[width=12cm]{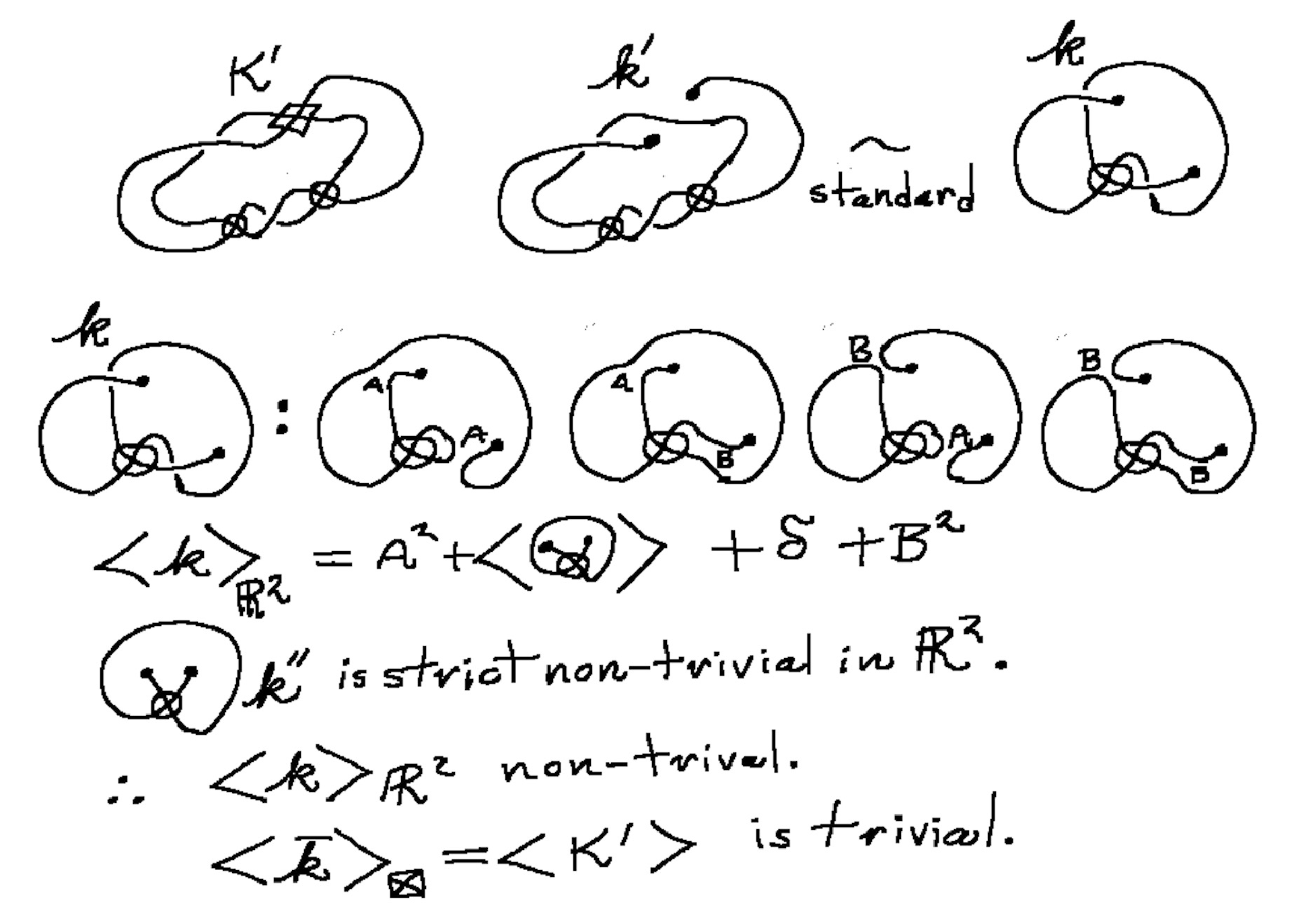}
     \end{tabular}
     \caption{\bf Strict Non-Triviality and Closure Triviality}
     \label{strictnon}
\end{center}
\end{figure}

Consider Figure~\ref{strictnon}. In the upper left of the this figure we show the multi-virtual knot $K'$ that is analyzed in Figure~\ref{virt} and Figure~\ref{nontriv}. In those figures we have shown that $K'$ has trivial generalized bracket polynomial but that nevertheless $K'$ is non-classical and non-trivial via a computation of a generalized Alexander polynomial in Figure~\ref{nontriv}. We now use this information to examine the strict planar knotoid $k'$ in Figure~\ref{strictnon}. First of all, note that $K'$ is the box closure of $k'$ and so the non-triviality of $K'$ proves that $k'$ is non-trivial as a strict knotoid on $S^2$ and hence as a planar knotoid as well.\\

Figure ~\ref{strictnon} illustrates a standard knotoid equivalence of $k'$ with the knotoid $k$ that is clearly trivial as an $S^2$ strict knotoid since a classical Reidemeister two move will undo it on the two-sphere. Nevertheless, $k$ is non-trivial as a strict planar knotoid. The figure illustrates the proof of this fact by performing a strict bracket expansion wherein we find that $\langle k \rangle = A^2 + B^2 + \delta + \langle k'' \rangle = \langle k'' \rangle$ (since $\delta = -A^2 - B^2$ with $B = A^{-1}$). The very simple planar knotoid $k''$ consists in two neighboring endpoints crossed virtually and surrounded by their own arc. $k''$ is a non-trivial planar knotoid in the strict category exactly because no detour move is available to simplify it. We omit the details of the proof of this final fact. This example shows how non-trivial entities without crossings that are different from the Turaev variables $U_{n}$ now begin to appear in the virtual strict context.\\

\begin{figure}
     \begin{center}
     \begin{tabular}{c}
     \includegraphics[width=12cm]{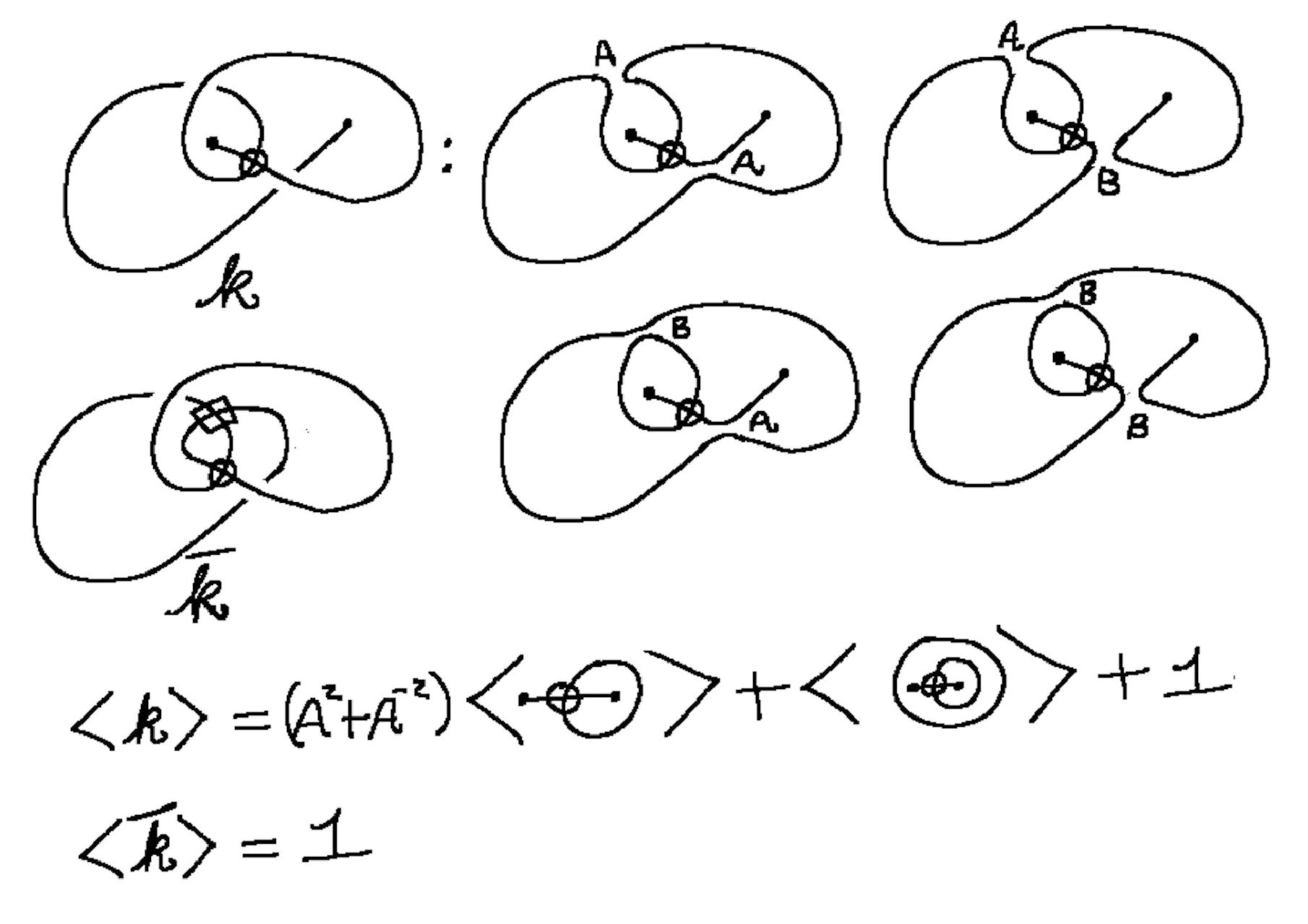}
     \end{tabular}
     \caption{\bf Generalized Loop Bracket}
     \label{genloop}
\end{center}
\end{figure}

Figure~\ref{genloop} shows a knotoid named $k$ in that figure (different from the previous figure) and its bracket expansion in the strict context. Note that there are states consisting in arcs intersected virtually by loops as well as encircled by loops. Thus we have irreducible states that are again different from the nested loops around arcs. These new states derive directly from strict irreducibility. Note that the box closure of $k$, denoted $\bar{k}$ in the figure, does not detect the non-triviality of $k$ as it is a trvial mutli-knotoid.\\

\begin{figure}
     \begin{center}
     \begin{tabular}{c}
     \includegraphics[width=12cm]{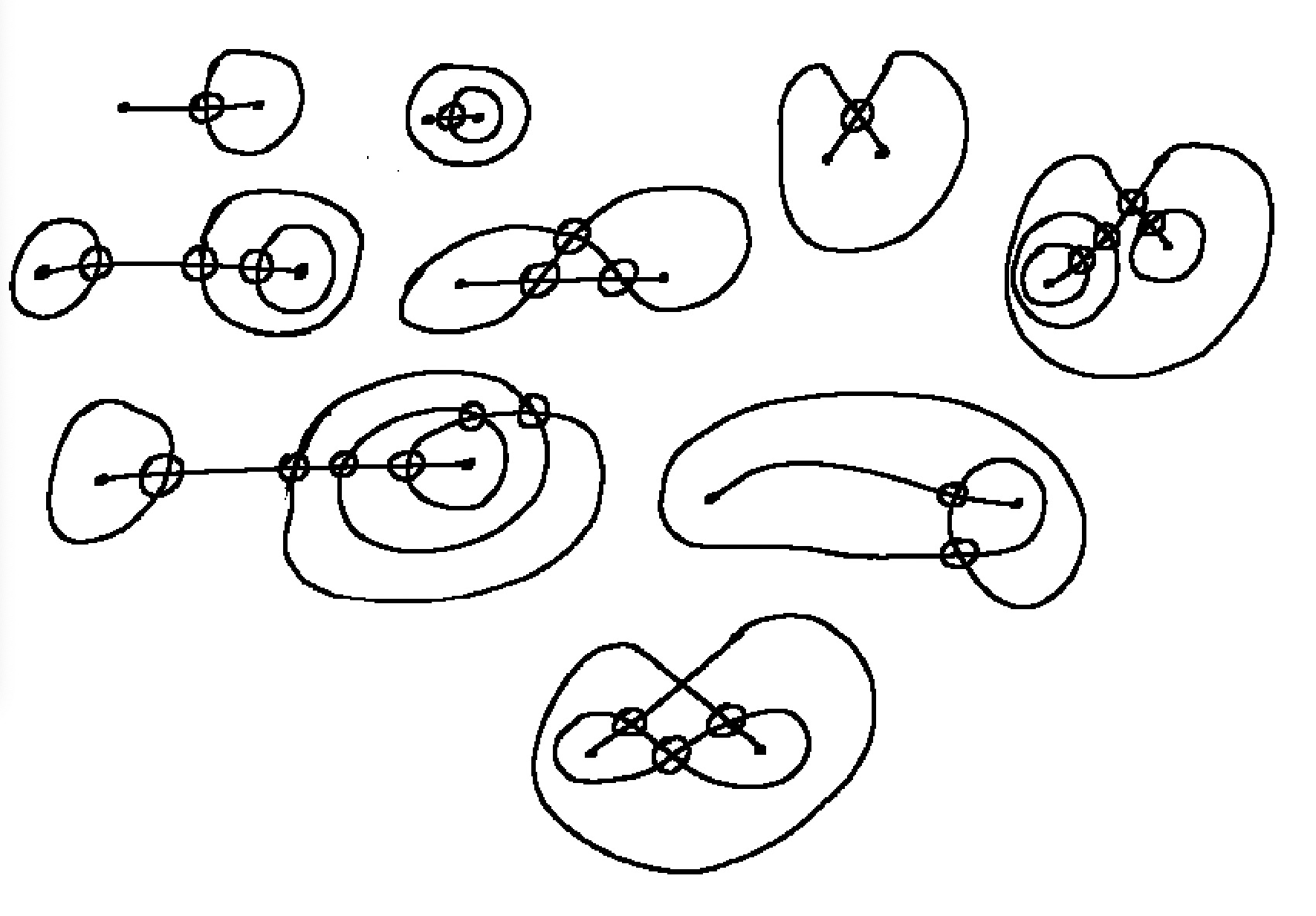}
     \end{tabular}
     \caption{\bf Examples of Planar Reduced States}
     \label{examples}
\end{center}
\end{figure}

In Figure~\ref{examples} we illustrate examples of reduced planar states. Such states need to be classified up to strict planar virtual isotopy. This figure illustrates the prolixity of such states even in the case of one virtual crossing type.
These examples are the beginning of an exploration of planar multi-virtual linkoids and knotoids that will be continued in subsequent work.\\

\section{Strict Virtual Polar Linkoids}

\begin{figure}
     \begin{center}
     \begin{tabular}{c}
     \includegraphics[width=12cm]{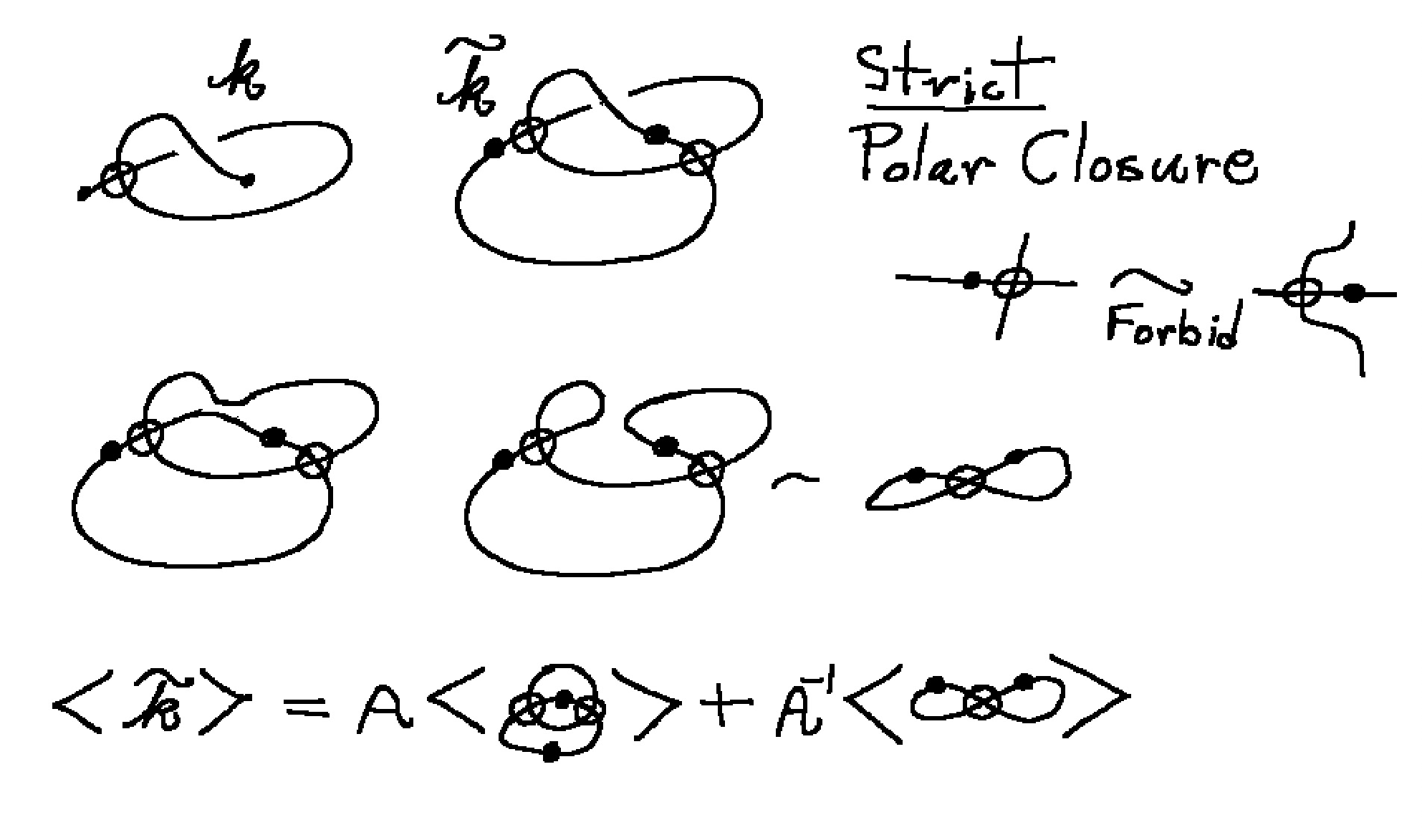}
     \end{tabular}
     \caption{\bf Strict Polar Linkoids and Closures}
     \label{strictpolar}
\end{center}
\end{figure}

\begin{figure}
     \begin{center}
     \begin{tabular}{c}
     \includegraphics[width=12cm]{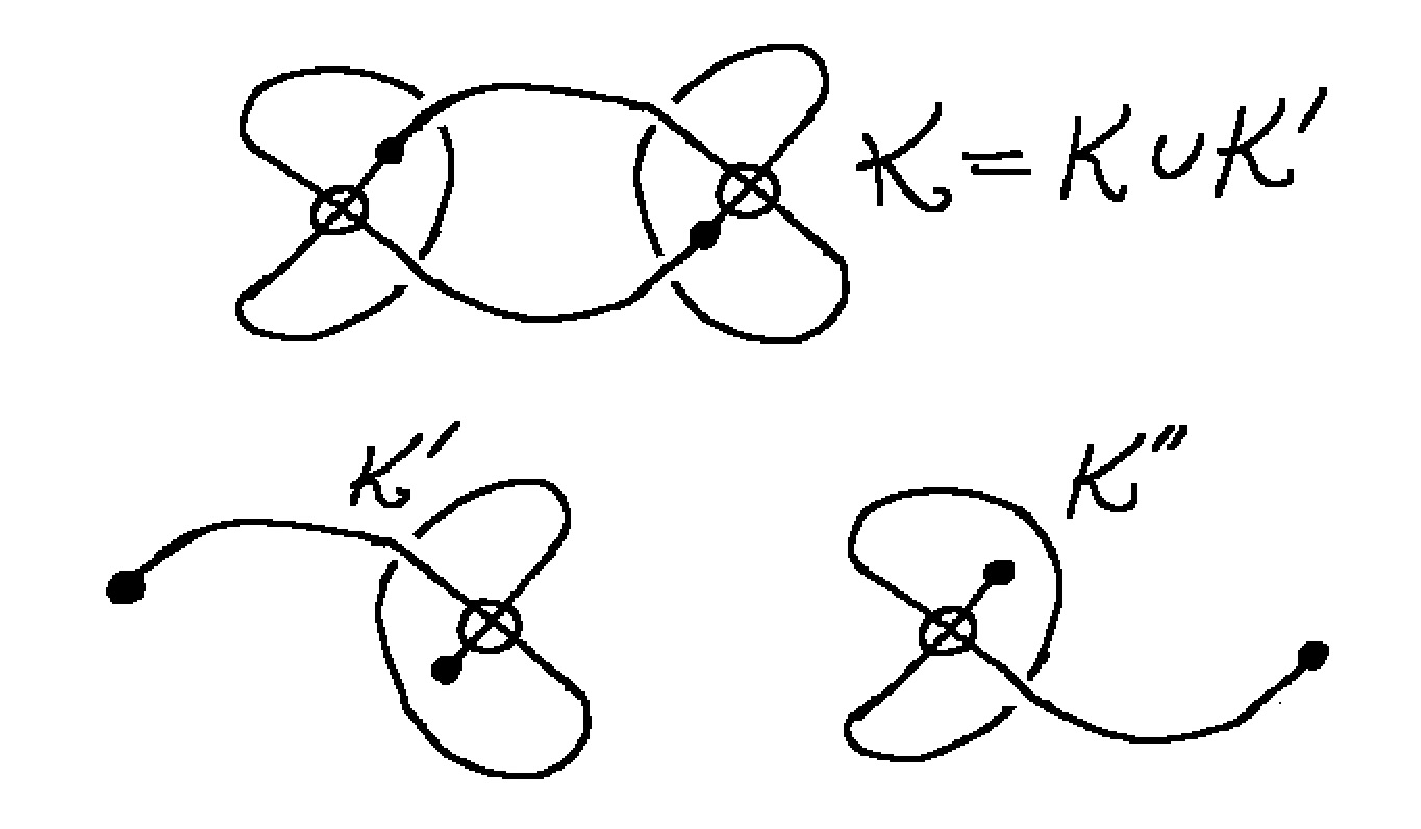}
     \end{tabular}
     \caption{\bf A Virtual Polar Linkoid and its Linkoid Parts}
     \label{polarparts}
\end{center}
\end{figure}

\begin{figure}
     \begin{center}
     \begin{tabular}{c}
     \includegraphics[width=12cm]{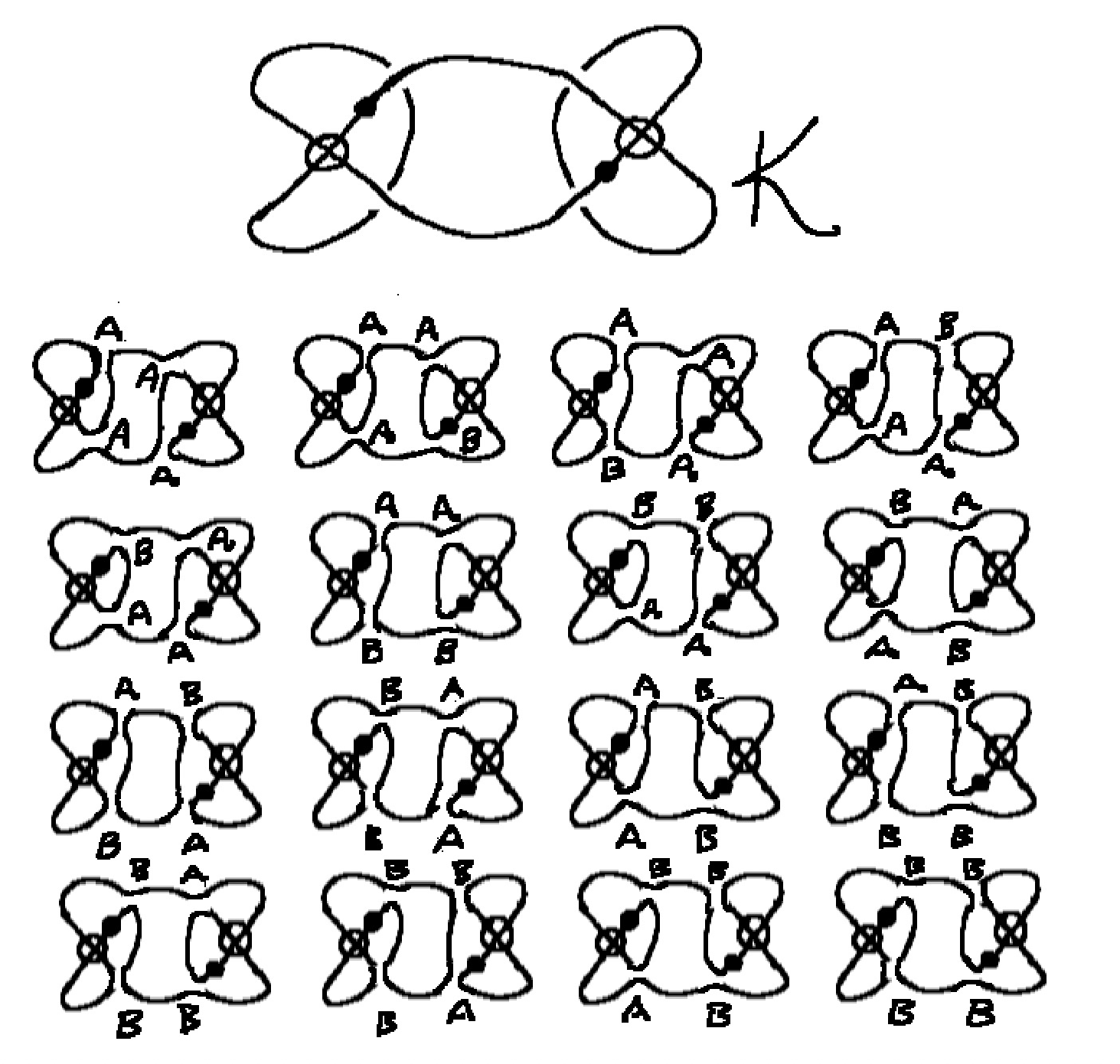}
     \end{tabular}
     \caption{\bf Bracket States for A Virtual Polar Linkoid}
     \label{polarstates}
\end{center}
\end{figure}

\begin{figure}
     \begin{center}
     \begin{tabular}{c}
     \includegraphics[width=12cm]{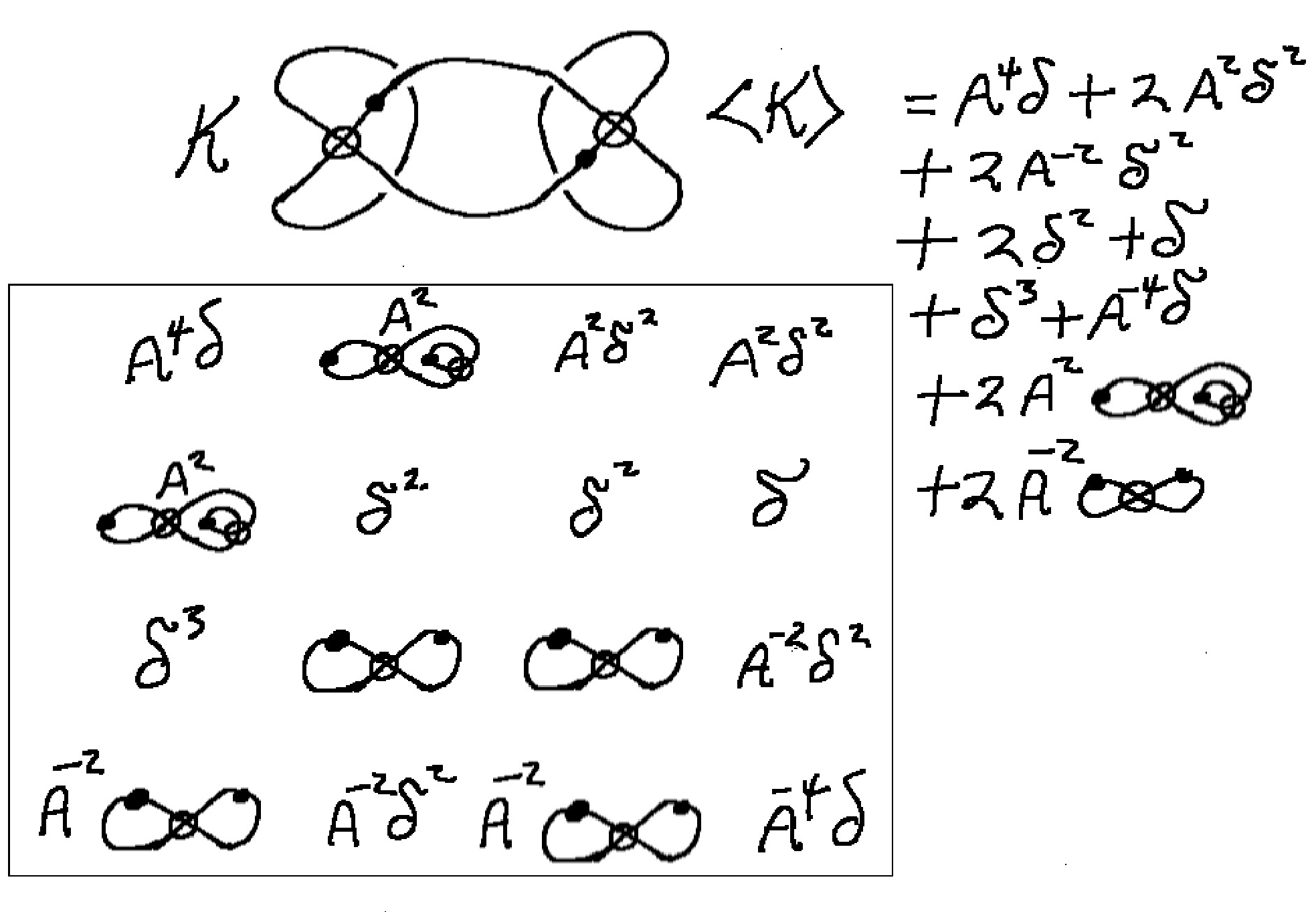}
     \end{tabular}
     \caption{\bf Reduced States and Bracket Calculation}
     \label{polarbracket}
\end{center}
\end{figure}

In this section we consider a different sort of closure for a knotoid or linkoid. in this kind of virtual closure we retain the endpoints of the knotoid in the virtual link that is the closure, and we, at the outset, do not introduce new 
virtual crossing types. By retaining the endpoints we have a virtual link with degree two nodes in some of its arcs. We follow \cite{Bataineh} in calling such virtual links {\it polar virtual links}. In \cite{Bataineh} the polar links have only classical crossings, so the virtual polar links are an innovation of this paper. In the original polar links and in virtual polar links, it is not allowed to move arcs across the degree two nodes. We will use strict equivalence on the virtual polar links so that all moves, including virtual moves, are performed in local move fashion as in strict knotoid equivalence. Thus we will work with strict virtual polar links. In the examples considered here we have single virtual crossings in our virtual polar links, but the theory accomodates multi-virtual crossings in the same way. \\

Figure~\ref{strictpolar} illustrates the polar closure of a strict knotoid $k$ and the calculation of the bracket polynomial of the corresponding strict virtual polar knot. Note that the virtual polar knot $\tilde{k}$ has only round virtual crossings.
Nevertheless, the bracket polynomial of $\tilde{k}$ is non-trivial becuase its states contain degree two nodes that are obstructions to virtual detour simplications that do not move over these nodes. It is not hard to prove that the two states illustrated in this figure are irreducible under strict detour moves (by puncturing the plane at the degree two nodes and considering the homotopy classes of the remaining loops in the complement of these nodes).\\

In Figure~\ref{polarparts} we show a strict virtual polar knot $K$ that is the union of two strict virtual knotoids $K'$ and $K''.$ These two knotoids are equivalent to one another in their own form and we have worked with them in the previous section of the paper, proving that they are non-trivial strict knotoids. Since equivalence of strict polar knots implies the equivalence of any subknotoids in them, we know that $K$ is non-trivial because $K'$ and $K''$ are non-trivial.\\

Another way to prove that $K$ is non-trivial is shown in Figure~\ref{polarstates} and Figure~\ref{polarbracket} where, in the first figure we show the full 16-state expansion for the bracket of $K,$ and in the second figure we simplify these states and collect them into the generalized bracket of $K.$ The resulting generalized bracket has irreducible states the show the strict non-triviality of $K.$ This method can be compared with the method of the previous paragraph where we examined sub-knotoids of the virtual polar knotoid.\\

This section has introduced the concept of  strict virtual polar knots and links and how their properties can be understood in terms of the theory of strict linkoids. In turn this paper has shown how strict linkoids can be studied in terms of the theory of multi-virtual
knots, links and their corresponding knotoids and linkoids. Much more work can be done in these domains and this paper is intended as an introduction to these investigations.\\

\section{Problems and Prospects}
In this section we summarize problems stated in the body of the paper and ask further questions.\\

\begin{enumerate}
\item In section 2.4 and Figure~\ref{Handles} we discuss how multiple virtual diagrams can be regarded as diagrams of links drawn on a sphere with attached handles, one handle for each virtual crossing. The handles carry labels corresponding labeled virtual crossings. Handle stabilization is restricted to interactions with handles of the same label type. This point of view deserves further exploration. There are variations that will be of interest. For example, if $K$ is a diagram with two virtual types, call them ``round" and ``square", then we could add handles only at the round virtual crossings, obtaining a single-type virtual diagram on a surface whose genus is  determined by the configuration of round crossings. Single-type virtual diagrams on surfaces of arbitrary genus  and indeed multi-virtual diagrams on arbitrary surfaces are natural generalizations of virtual knot theory that  deserve further exploration. Many invariants immediately generalize to this context and can be compared with our known invariants of multi-virtuals in the two sphere and the plane. Much more can be done in this domain.
\item  In section 2.5 we define the generalized bracket invariant for multi-virtuals.
The generalized bracket has state expansion 
$$<K> = \Sigma_{S} A^{a(S)}B^{b(S)} \langle S \rangle$$ 
where $S$ is a state of the diagram $K$ obtained by choosing a smoothing for each 
crossing of $K$ and $a(S)$ is the number of smoothings of type $A,$ while $b(S)$ is the number of smoothings of type $B.$ Smoothing types correspond to the convention shown in 
Figure~\ref{EFF2} where the local regions of a crossing are labeled $A$ and $B$ so that the $A$ regions are swept by a counterclockwise turn of the overcrossing arc. A smoothing is of type $A$ if the two $A$ regions are joined by the smoothing. Each state $S$ in the virtual bracket expansion is a diagram with multi-virtual crossings. For the generalized bracket, we take $\langle S \rangle$ to be the multi-virtual class of the state $S$ with the caveat that any free disjoint circle in $S$ (or in a diagram equivalent to $S$) is evaluated as $\delta.$ Thus we have that $\langle S O \rangle = \delta \langle S \rangle.$ In \cite{MVKT} we consider closed loop examples of this problem and indeed Figure~\ref{onecomp} illustrates a conjecture that a single-loop state with exactly two distinct virtual crossings can be undone by virtual moves, while Figure~\ref{EFF8} shows how a two loop state with two distinct virtual crossings is not equivalent to two disjoint loops, by using the chromatic bracket. The rest of the paper uses variations of the bracket so that $\langle S \rangle$ depends on the different contexts and equivalence relations. A good case in point is shown
in Figure~\ref{examples} where these diagrams are to be classified up to strict virtual equivalence for planar linkoids. This means that only compositions of planar local virtual moves are allowed, and no move can shift an arc across an endpoint. The generalized bracket is defined for such linkoids, but we do not yet have a complete classification of the states $S.$
\item Many problems arising in the present paper suggest that the subject of diagrams without classical crossings and only multiple types of virtual crossings (taken either strictly or with full detours) is a subject in its own right. Classify such diagrams in the plane and on other surfaces.
\item Figure~\ref{virt} and Figure~\ref{nontriv} recall the virtualization method and how it can be used to produce multi-virtual knot diagrams and strict 1-linkoid diagrams with unit Jones polynomial. In the case of a single virtual type we have shown that no non-trivial classical knot diagram with unit Jones polynomial can be produced by this method \cite{DKK}. Can this result be extended to the multi-virtual case? Is there a non-trivial extension of Khovanov Homology for virtual knots and links to multi-virtual knots, links and linkoids?
\item The link $\bar{k}$ in Figure~\ref{invisible} is undetectable by the generalized bracket polynomial for multi-virtual links. Find more examples of this kind with an eye toward understanding the corresponding phenomenon of invisibility with respect to the Jones polynomial for classical links.
\item In this paper we have generalized polar knots and links \cite{Bataineh} to multi-virtual polar knots and links, in the process showing how strict multi-virtual equivalence is a natural way to study sub-linkoids of polar links and knots. The idea is that once one constrains the isotopy of a given closed diagram by declaring that one cannot move arcs across certain degree two nodes, then the strict isotopy classes of sub-linkoids of that diagram  will give invariant information. More work can be done here. In particular one can define more invariants of multi-virtual polars, including new state sum invariants and including generalizations of quandle invariants. This will be the subject of subsequent work.

\end{enumerate}

\end{document}